\documentclass[a4paper,11pt,bibtotoc]{scrartcl}

\usepackage[utf8]{inputenc}

\usepackage{pst-node}
\usepackage{tikz-cd} 
\usepackage{wrapfig}
\usepackage{caption}

\usepackage{amsthm}

\usepackage{amsmath,amssymb}

\usepackage{dsfont}

\usepackage{mathtools}

\usepackage{mathrsfs}

\usepackage[new]{old-arrows}

 \usepackage{comment}

\usepackage[
backend=biber,
style=alphabetic,
sorting=nyt,
backref=true
]{biblatex}

\DefineBibliographyStrings{english}{%
  backrefpage = {p.},
  backrefpages = {pp.},
}

\usepackage{todonotes}

 \usepackage {hyperref}
 \hypersetup{
    colorlinks=true,
    linkcolor=blue,
    filecolor=magenta,      
    urlcolor=cyan,
    citecolor=magenta,
    pdftitle={Overleaf Example},
    pdfpagemode=FullScreen,
    }

\usepackage{MnSymbol}

\usepackage{placeins} 

\usepackage{etoolbox} 

\usepackage[legalpaper, a4paper,
margin=1in]{geometry}

\robustify{\underset}

\DeclareRobustCommand{\bbone}{\text{\usefont{U}{bbold}{m}{n}1}}
\DeclareMathOperator{\EX}{\mathbb{E}} 
\DeclareMathOperator{\PR}{\mathbb{P}} 
\newcommand{\defeq}{\overset{\mathrm{def}}{=\joinrel=}}

\DeclarePairedDelimiter{\floor}{\lfloor}{\rfloor}
\DeclarePairedDelimiterX{\inner}[1]{\langle}{\rangle}{#1}
\def\acts{\curvearrowright} 
\def\Bcovers{\le_{\overset{\twoheadrightarrow}{B}}}
\def\leff{\overset{*}{\le}}
\def\lealg{\le_{\textup{alg}}}

\def\multiset#1#2{\ensuremath{\left(\kern-.3em\left(\genfrac{}{}{0pt}{}{#1}{#2}\right)\kern-.3em\right)}}

\newcommand{\Bconv}{\underset{B}{*}}
\newcommand{\convalg}{\underset{\textup{alg}}{*}}

\DeclareMathOperator{\C}{\mathbb{C}}

\DeclareMathOperator{\F}{\mathbb{F}}

\DeclareMathOperator{\N}{\mathbb{N}}
\DeclareMathOperator{\Q}{\mathbb{Q}}
\DeclareMathOperator{\R}{\mathbb{R}}
\DeclareMathOperator{\Ss}{\mathbb{S}}

\DeclareMathOperator{\Z}{\mathbb{Z}}

\newtheorem{theorem}{Theorem}[section]
\newtheorem{corollary}[theorem]{Corollary}

\newtheorem{lemma}[theorem]{Lemma}

\newtheorem{proposition}[theorem]{Proposition}
\newtheorem{definition}[theorem]{Definition}
\newtheorem{example}[theorem]{Example}
\newtheorem{remark}[theorem]{Remark}
\newtheorem{conjecture}[theorem]{Conjecture}
\newtheorem{observation}[theorem]{Observation}

\title{Word Measures on Wreath Products I}
\author{Yotam Shomroni }
\date{April 2023}

\begin{document}

\maketitle

\begin{abstract}
    Every word $w$ in the free group $F_r$ of rank $r$ induces a probability measure (the $w$-measure) on every compact group $G$, by substitution of Haar-random $G$-elements in the letters. 
    This measure is determined by its Fourier coefficients: the $w$-expectations $\EX_w[\chi]$ of the irreducible characters of $G$.
    For every compact group $G$, the wreath product with the symmetric group $G\wr S_n$ has some natural irreducible characters $\chi$, and we approximate $\EX_w[\chi]$ for every word $w\in F_r$, revealing new automorphism-invariant quantities of words that generalize the primitivity rank $\pi(w)$.
    This generalizes previous works by Parzanchevsky-Puder and Magee-Puder.
    We demonstrate applications to automorphism groups of trees, investigate properties of the new invariants, and show polynomial decay of $\EX_w[\chi]$ also for wreath products with more general actions.
\end{abstract}

\tableofcontents 

\section{Introduction}
\label{section_intro}

In this paper we study word measures on wreath products of compact groups with the symmetric group $S_n$.
Primarily, we bound the $w$-expectations of some irreducible characters. 

We start by explaining the notions of wreath product and word measures.
The \textbf{wreath product} of a group $G$ with the symmetric group $S_n$ is $G\wr S_n\defeq G\wr_{[n]} S_n \defeq G^n\rtimes S_n$, where $S_n$ acts on $G^n$ by permuting the indices $[n]\defeq \{1, \ldots, n\}$, i.e.\ $\sigma.(g_1, \ldots, g_n) = (g_{\sigma(1)}, \ldots, g_{\sigma(n)})$. 
The elements are $\left\{(v, \sigma): v\in G^n, \, \sigma\in S_n\right\}$ and the product is $\left(v_1, \sigma_1\right) \cdot \left(v_2, \sigma_2\right) = \left(v_1 \cdot (\sigma_1.v_2), \sigma_1\cdot \sigma_2\right)$. 
An element can also be thought of as a monomial matrix (every row and column has a unique non-zero entry) with $G$-elements, e.g.\
\begin{equation*}
    \begin{pmatrix}
    0 & i & 0 \\
    0 & 0 & 1 \\
    -i & 0 & 0
    \end{pmatrix} \in \{\pm 1, \pm i\} \wr S_3, \quad
    \begin{pmatrix}
    g_1 & 0 & 0 \\
    0 & 0 & g_2 \\
    0 & g_3 & 0
    \end{pmatrix} \in G\wr S_3 \quad (g_i\in G).
\end{equation*}

\subsection*{What are word measures?}

Given a word $w\in F_r$ in a free group $F_r = \textup{Free}\left(\{b_1, \ldots, b_r\}\right)$, and a group $G$, we get a map (not necessarily a homomorphism)\footnote{Even though not being homomorphisms themselves, word maps commute with homomorphisms. In the language of categories, word maps are precisely the natural transformations $\textbf{forg}^r\to \textbf{forg}$, where $\textbf{forg}\colon \textbf{Grp}\to \textbf{Set}$ is the forgetful functor and $\textbf{forg}^r(G) = G^r$ gives the $r$-fold Cartesian product, and similarly $\textbf{forg}^r(f) = (f, \ldots, f)$ for every homomorphism $f$. } 
$w\colon G^r\to G$, called a word map. 
For example, $b_1 b_2^2 b_1^{-1} b_2^{-1}$ maps $(g, h)\in G^2 \mapsto gh^2 g^{-1}h^{-1}\in G$.
If $G$ is compact, it has an invariant probability measure $\mu$ (the Haar measure), and similarly $\mu^{\times r}$ is an invariant probability measure on $G^r$. The pushforward measure $\mu_w \defeq w_*(\mu^{\times r})$ on $G$ is called the $w$-measure on $G$.
Equivalently, by the universal property of the free group, homomorphisms $F_r\to G$ correspond to functions $\{b_1, \ldots, b_r\}\to G$, which can be identified with $G^r$, so $\textup{Hom}(F_r, G)$ inherits the unique shift-invariant probability measure $U$ from $G^r$. 
Given a random homomorphism $\alpha\sim U(\textup{Hom}(F_r, G))$, the image $\alpha(w)$ distributes according to the $w$-measure. If $G$ is finite, this reduces to 
\[\forall g\in G:\,\, \mu_w(g) = \frac{1}{|G|^r} \left|\left\{(g_1, \ldots, g_r)\in G^r: w(g_1, \ldots, g_r) = g\right\}\right|. \]
For example, the word $b_1\in F_r$ induces the Haar (uniform) probability measure on any compact group.
Given a function $f\colon G\to \C$, we wish to compute its expectation according to the $w$-measure:
\[ \EX_w[f] \defeq \EX_{\alpha}[f(\alpha(w))] = \int_{g\in G} f(g) \mathrm{d}\mu_w(g). \]

This expectation captures interesting information about both $w$ and $f$.
As a toy example, consider $G = C_m$, the cyclic group of order $m$, and $f\colon C_m \hookrightarrow \Ss^1$ the standard embedding $f(x) = e ^ {\frac{2\pi i x}{m}}$.
Given a basis $B = \{b_1, \ldots, b_r\}$ of the free group $F_r$, define $\nu_i\in \textup{Hom}(F_r, \Z)$ for every $i\in [r]$ by letting $\nu_i(b_j)\defeq \bbone_{i=j}$ for every $j\in [r]$. 
(This determines $\nu_i$ uniquely).
Then 
\begin{equation}
\label{equation_expectation_of_cyclic_embedding}
\EX_w[f] =
    \begin{cases}
    1 & \textrm{ if } \nu_i(w)\equiv 0 \textup{ (mod }m\textup{)} \textrm{ for every }i,\\
    0 & \textrm{otherwise.}
    \end{cases}
\end{equation}
Indeed, by independence of the letters, and since the sum of a non-trivial subgroup of the $m^{th}$-roots of unity vanishes,
\[\EX_w[f] = \prod_{i=1}^r \EX_{b_i^{\nu_i(w)}}[f] = \prod_{i=1}^r \bbone_{\nu_i(w) = 0 \textup{ (mod }m\textup{)}}.\]


A much more interesting example arises when $G = S_n$ is the symmetric group, and $f = \#\textup{fix}$ is the natural character that counts fixed points of permutations.
In \cite[Theorem 1.8]{PP15}, Puder and Parzanchevsky proved an approximation theorem for $\EX_w[f]$:
\begin{equation*}
    \label{equation_PP15_thm_main_result}
    \EX_w[f] = 1 + |\textup{Crit}(w)|n^{1-\pi(w)} + O\left(n^{-\pi(w)}\right),
\end{equation*}
where $\pi(w), \textup{Crit}(w)$ are the \textbf{primitivity\footnote{{By Nielsen-Schreier theorem, a subgroup of a free group is free. 
A primitive element in $H$ is a part of a basis of $H$, or equivalently, in the $\textup{Aut}(H)$-orbit of $h_1$ for some basis $\{h_1, \ldots, h_k\}$ of $H$.}} rank} of $w$ and the \textbf{critical subgroups}, defined as follows:
\begin{definition}
\label{def_witnesses}
\label{def_primitivity_rank}
(Primitivity Rank; \cite[Definition 1.7]{Puder_2014}) 
Let $w\in F_r$ be a word in a free group.
The \textbf{primitivity rank} of $w$ is
\[\pi(w) \defeq \min\left\{\textup{rk}(H): w\in H\le F_r, \, w\textrm{ is a non-primitive element in }H\right\},\]
with the convention $\min\emptyset = \infty$.
The \textbf{critical subgroups} are the subgroups achieving the minimum:
\[\textup{Crit}(w) \defeq \left\{H\le F_r: w\in H\textrm{ is a non-primitive element in }H, \, \textup{rk}(H)=\pi(w)\right\}.\]
\end{definition}

The image of $\pi$ is $\{0, 1, \ldots, r, \infty\}$. The following table gives examples for $\pi(w), \textup{Crit}(w)$:

\begin{table}[ht!]
\centering
\begin{tabular}{||c | c | c ||} 
 \hline
 Description of $w$ 
 & $\pi(w)$ 
 & $\textup{Crit}(w)$ \\[0.7ex] 
 \hline\hline
 $w=1$ & 0 & $\{\inner*{1}\}$ \\[0.7ex] 
 $w$ is a proper power & 1 & $\{\inner*{u}: \inner*{w} \lvertneqq \inner*{u}\}$  \\[0.7ex] 
 $[b_1, b_2]$ & 2 & $\{\inner*{b_1, b_2}\}$ \\[0.7ex] 
 $b_1^2 b_2^2$ & 2 & $\{\inner*{b_1, b_2}\}$ \\[0.7ex] 
 $\vdots$ & $\vdots$ & $\vdots$ \\[0.7ex] 
 $b_1^2 \ldots b_k^2$ & k & $\{\inner*{b_1, \ldots, b_k}\}$ \\[0.7ex] 
 $w$ is primitive & $\infty$ & $\emptyset$ \\ [1ex] 
 \hline
\end{tabular}
\caption{Primitivity Rank and Critical Subgroups}
\label{table:primitivity_rank}
\end{table}
\FloatBarrier

In \cite{magee2021surface}, Magee and Puder studied word measures on wreath products of the form $C_m\wr S_n$ for $m\in \{2, 3, \ldots\}\sqcup \{\infty\}$, where $C_{\infty} = \Ss^1$ is defined to be the unit circle group (and recall the notation $C_m = \Z/m$). By representing elements of $C_m\wr S_n$ as monomial matrices (as $(v, \sigma)\mapsto (v_i\cdot \bbone_{\sigma(i) = j})_{i, j = 1}^n$), we get a natural embedding $C_m\wr S_n\hookrightarrow \textup{GL}_n(\C)$. This is a "standard" $n$-dimensional irreducible representation, whose character is just the trace of the monomial matrix, similarly to the fixed-point counting for $m = 1$.
To describe the result of \cite{magee2021surface} we need some more definitions. 

\begin{definition}
\label{def_K_m(H)}
Let $H$ be a free group and let $m\in \N$. Define
\[ K_m(H)\defeq \ker\left(H\twoheadrightarrow H^{ab}/m \right)\]
(where $H^{ab}/m\defeq H/\llangle [H, H], H^m \rrangle \cong C_m^{\textup{rk}(H)}$).
This is the characteristic subgroup of $H$ obtained by taking only words in $H$ that use each letter (for some basis of $H$) $m\Z$ times in total\footnote{Equivalently, $w\in K_m(H)$ if and only if the presentation complex of $H/\inner*{\inner*{w}}$ has non-trivial second homology with coefficients in $\Z/m$.}.
Similarly, define $K_{\infty}(H)\defeq [H, H]$.
\end{definition}

\begin{example}
In the group $F_3 = \textup{Free}(\{a, b, c\})$, the word $w = a^2 b c b c^{-1}$ is sent to $w^{\textup{ab}} = (2, 2, 0)$ under the abelianization map to $F_3^{\textup{ab}} \cong \Z^3$, so $w\in K_2(F_3)$. On the other hand, $w\not \in K_2(\inner{w})$.
\end{example}

\begin{definition}
\label{def_pi_Cm}
(\cite[Definition 1.10]{magee2021surface})
For every $m\in \N_{\ge 2}$, define $\pi_{C_m}\colon F_r\to \{0, 1, \ldots, r, \infty\}$ (see Remark~\ref{remark_composition_with_balanced_morphisms} for an explanation about the range of $\pi_{\phi}$) by
\[ \pi_{C_m}(w) \defeq \min\{ \textup{rk}(H): H\le F_r, w\in K_m(H) \}. \]
Also define $\textup{Crit}_{C_m}(w)\defeq \{H\le F_r: w\in K_m(H),\,\, \textup{rk}(H) = \pi_{C_m}(w)\}$ as the set of subgroups achieving the minimum.
For $m = \infty$ define $\pi_{\Ss^1}(w) = \min\{\textup{rk}(H): H\le F_r, w\in [H, H]\}$.
\end{definition}

Now for every $m\in \{2, 3, \ldots\}\sqcup \{\infty\}$, \cite[Theorem 1.11]{magee2021surface} approximates the expected trace of a $w$-random element of $C_m\wr S_n$: 
\[ \EX_w[\textup{tr}] = |\textup{Crit}_{C_m}(w)|\cdot n^{1-\pi_{C_m}(w)} + O\left(n^{-\pi_{C_m}(w)}\right). \]
Note that if $d\divides m$ then $K_m(H)\le K_d(H)$ so $\pi_{C_d}(w)\le \pi_{C_m}(w) \le \pi_{\Ss^1}(w)$.
Also note that $w\in K_m(H)$ implies (when $m\ge 2$) that $w$ is not primitive in $H$, so $\pi(w)\le \pi_{C_m}(w)$.
In this sense, we may regard such $H$ as a witness for the imprimitivity of $w$.

\subsection*{Statement of the Main Result}

Now let $G$ be again an arbitrary compact group, and recall the notation $U$ for the Haar (uniform) measure on compact groups.
Inspired by \cite{magee2021surface}, for every word $w$ we can investigate many more "imprimitivity detecting subgroups", generalizing groups $H$ such that $w\in K_m(H)$.
For a class function $f\colon G\to \C$ and $w\in H\le F_r$ (here $H$ is a finitely generated subgroup), denote 
$\EX_{w\to H}[f]\defeq \EX_{\alpha\sim U(\textup{Hom}(H, G))}[f(\alpha(w))]$.
Recall that for every compact group $G$ we have the standard normalized inner product on $\C^{G}$, that is $\inner*{\phi, \psi}\defeq \int_{g\in G} \phi(g) \overline{\psi}(g)$ (integrating with respect to the Haar measure), and that the irreducible characters form an orthogonal basis for class functions\footnote{By Peter–Weyl theorem, the characters of the irreducible representations of $G$ form a Hilbert basis for the space of square-integrable class functions on $G$.}.
So in other words, if $\mu_{w\to H}$ is the $w$-measure on $G$ when we think of $w$ as an element of $H$, then 
$\EX_{w\to H}[f]\defeq \int_{g\in G} f(g) \textup{d}\mu_{w\to H}(g)$.
Note that if $w$ is primitive in $H$ then $\EX_{w\to H}[f] = \inner*{f, \textbf{1}} = \int_{g\in G}f(g)$.
(We denote by $\textbf{1}\colon G\to \{1\}$ the trivial character, which is the constant $1$).
We denote the set of irreducible characters of a group $G$ by $\hat{G}$, the set of all characters by $\textup{char}(G)$, and the set of conjugacy classes by $\textup{conj}(G)$.

To state the main result, we start by generalizing the primitivity ranks $\pi(w)$ and $\pi_{C_m}(w)$. 
Since for every primitive word $w\in F_r$ and a non-trivial irreducible character $\phi\in \hat{G}-\{\textbf{1}\}$ of a compact group $G$ we have $\EX_{w\to H}[\phi] = 0$, for every non-primitive word $w$ we may regard subgroups $H$ that violate this condition as witnesses for the imprimitivity of $w$.
\begin{definition}
\label{def_witnesses_no_alg}
($\phi$-Witnesses of Imprimitivity)
Let $G$ be a compact group and $F_r$ the free group of rank $r$.
For every word $w\in F_r$ and an irreducible character $\phi\in \hat{G}$, we define the $\phi$-primitivity rank of $w$ as
\begin{equation}
    \label{eq_phi_prim_alternative}
        \pi_{\phi}(w) \defeq 
        \begin{cases}
        \min\{\textup{rk}(H): w\in H \le F_r, \,\EX_{w\to H}[\phi] \neq 0\} & \textup{ if } \phi\neq \textbf{1},\\
        \pi(w) & \textup{ if } \phi = \textbf{1}.
        \end{cases}
    \end{equation}
We also define the $\phi$-critical groups of $w$, which are the subgroups that achieve the minimum:
    \begin{equation*}
        \textup{Crit}_{\phi}(w) \defeq \begin{cases}
        \{H\le F_r: w\in H, \,\EX_{w\to H}[\phi] \neq 0, \textup{rk}(H) = \pi_{\phi}(w)\} & \textup{ if } \phi\neq \textbf{1},\\
        \textup{Crit}(w) & \textup{ if } \phi = \textbf{1}.
        \end{cases}
    \end{equation*}
Finally, we define the $\phi$-critical value of $w$ as
    \[\mathscr{C}_{\phi}(w) \defeq \sum_{H\in \textup{Crit}_{\phi}(w)} \EX_{w\to H}[\phi]. \]
\end{definition}

The definition $\pi_{\textbf{1}} = \pi$ looks ad hoc, but in Proposition~\ref{prop_equiv_def_pi_phi} we give the "correct" equivalent definition of $\pi_{\phi}$ which captures both cases at once, using the concept of algebraic extensions of free subgroups.

When $G$ is finite, $\mu_{w\to H}\colon G\to [0, 1]$ is a class function, so 
$\mu_{w\to H} = \sum_{\phi\in \hat{G}} \inner{\mu_{w\to H}, \phi} \phi = \frac{1}{|G|} \sum_{\phi\in \hat{G}} \EX_{w\to H}[\phi]\cdot \phi$.
Hence $\mu_{w\to H}$ is "built up" from all $\phi\in \hat{G}$ satisfying $\EX_{w\to H}[\phi] \neq 0$.

This definition generalizes $\pi, \pi_{C_m}$ and $\textup{Crit}, \textup{Crit}_{C_m}$, by the example in equation~\eqref{equation_expectation_of_cyclic_embedding}.
Indeed, let $\phi\colon C_m\to \C^{\times}$ be the inclusion map of the $m^{th}$ roots of unity (and $\Ss^1$ if $m=\infty$). Then for every word $w$, 
\begin{align*}
    & \pi_{\textbf{1}}(w) = \pi(w), && \textup{Crit}_{\textbf{1}}(w) = \textup{Crit}(w), && \mathscr{C}_{\textbf{1}}(w) = |\textup{Crit}(w)|, \\
    & \pi_{\phi}(w) = \pi_{C_m}(w), && \textup{Crit}_{\phi}(w) = \textup{Crit}_{C_m}(w), && \mathscr{C}_{\phi}(w) = \left|\textup{Crit}_{C_m}(w)\right|. 
\end{align*}
Note the automorphism-invariance of $\pi_{\phi}(w)$: if $\xi\in \textup{Aut}(F_r), \psi\in \textup{Aut}(G)$ then 
\[\pi_{\phi\circ \psi}(\xi(w)) = \pi_{\phi}(w),\quad \textup{Crit}_{\phi\circ \psi}(\xi(w)) = \textup{Crit}_{\phi}(w). \]
Indeed, for every $H$ containing $\xi(w)$, $\textup{rk}(H)=\textup{rk}(\xi^{-1}(H))$ and 
\[\EX_{\xi(w)\to H}[\phi\circ \psi] = \EX_{\alpha\sim U(\textup{Hom}(H, G))}[\phi(\psi \circ \alpha \circ \xi(w))] = \EX_{w\to \xi^{-1}(H)}[\phi] \] 
since pre-composition and post-composition with automorphisms do not change the Haar (uniform) distribution on homomorphisms (see for example \cite[Fact 2.5]{unitary_Magee_2019}).

\begin{example}
Let $G$ be a compact group and $\phi\in \hat{G}$ an irreducible character.
In \cite{frobenius1896gruppencharaktere}, Frobenius proved
$\EX_{[b_1,b_2]\to \inner*{b_1, b_2}}[\phi] = \frac{1}{\phi(1)}\neq 0$ 
and in \cite{frobeniusSchur1906}, Frobenius and Schur proved
\begin{equation*}
    \EX_{b_1^2\to \inner*{b_1}}[\phi] = 
    \begin{cases}
        1 & \textup{ if } \phi \textup{ is real, i.e.\ there is some } \Phi\in\textup{Hom}(G, \textup{GL}_d(\R)) \textup{ such that } \phi = \textup{tr}\circ \Phi,\\
        0 & \textup{ if } \phi \textup{ is complex, i.e.\ }\textup{Im}(\phi)\not\subseteq \R,\\
        -1 & \textup{ if } \phi \textup{ is quaternionic, i.e.\ otherwise}.
    \end{cases}
\end{equation*}
Moreover, $\textup{Crit}([b_1,b_2])=\{\inner*{b_1, b_2}\}$ and $\textup{Crit}(b_1^2) = \{\inner*{b_1}\}$.
Consequently $\textup{Crit}_{\phi}([b_1,b_2]) = \{\inner*{b_1, b_2}\}$ and $\textup{Crit}_{\phi}(b_1^2) \subseteq \{\inner*{b_1}\}$. 
Hence for $w\in \{b_1^2, [b_1, b_2]\}$ we have $\mathscr{C}_{\phi}(w) = \EX_w[\phi]$.
In particular, $\mathscr{C}_{\phi}(w)$ may obtain the values $\{0, 1, -1, \frac{1}{2}, \frac{1}{3}, \ldots\}$.
More examples appear in Section~\ref{section_witnesses}. 
\end{example}

For every compact group $G$ we wish to approximate $\EX_w[\chi]$, where $\chi$ is a "natural" character of the wreath product $G_n\defeq G\wr S_n$.
Recall the notation $[n] \defeq \{1, \ldots, n\}$.
As we handle actions other than $S_n\acts [n]$, we give a more general definition: 

\begin{definition}
\label{def_general_Ind_phi}
($\textup{Ind}_n\phi, \chi_{\phi, n}$)
Let $G$ be a compact group, and let $\Sigma_n\le S_n$ be a transitive permutation group\footnote{For example $\Sigma_n$ could be $S_n$, $\textup{GL}_N(\F_q)$ with $n=q^N-1$, $S_N$ with $n = \binom{N}{2}$, etc.} (that is $\Sigma_n\acts [n]$ transitively). 
Let $\phi\in \textup{char}(G)$. 
We give two equivalent definitions for $\textup{Ind}_n\phi\in \textup{char}\left(G\wr_{[n]} \Sigma_n\right)$ (we drop $\Sigma_n$ from the notation as it can always be inferred from the context): 
\begin{itemize}
    \item A formula: for every $v=(v(1), \ldots, v(n))\in G^n, \sigma\in \Sigma_n$:
    \[ \textup{Ind}_n\phi(v, \sigma) \defeq \sum_{i: \,\,\sigma(i)=i} \phi(v(i)). \]
    \item Character description: extend $\phi$ to $\Tilde{\phi}\colon G\wr_{[n]} \textup{stab}(n)\to \C$, where $\textup{stab}(n)\le \Sigma_n $ is the stabilizer of the point $n$, via $\Tilde{\phi}(v, \sigma) = \phi(v(n))$.
    Then 
    \[\textup{Ind}_n\phi \defeq \textup{Ind}_{G\wr_{[n]} \textup{stab}(n)}^{G\wr_{[n]} \Sigma_n } \Tilde{\phi}. \]
\end{itemize}
If $\phi$ is the character of $\Phi\colon G\to \textup{GL}_{\dim(\Phi)}(\C)$, the representation $\textup{Ind}_n\Phi$ that yields $\textup{Ind}_n\phi$ has the simple description
\begin{equation*}
    \begin{pmatrix}
    g_1 & 0 & 0 \\
    0 & 0 & g_2 \\
    0 & g_3 & 0
    \end{pmatrix} \in G\wr \Sigma_n \quad \mapsto
    \begin{pmatrix}
    \Phi(g_1) & 0 & 0 \\
    0 & 0 & \Phi(g_2) \\
    0 & \Phi(g_3) & 0
    \end{pmatrix} \in \textup{GL}_{n\cdot \dim(\Phi)}(\C).
\end{equation*}
When $\phi$ is irreducible, define for every $n\in \N$ the irreducible character $\chi_{\phi, n}\colon G\wr \Sigma_n \to \C$ by $\chi_{\phi, n} \defeq \textup{Ind}_n\phi- \bbone_{\phi = \textbf{1}}$.
\end{definition}

In Proposition~\ref{prop_Ind_phi_is_irred}, we will show that the two definitions coincide and give a well-defined character, and that if $\phi\in \hat{G}-\{\textbf{1}\}$ is non-trivial and irreducible, then so is $\chi_{\phi, n}$ (of dimension $\phi(1)n$).
When $\phi = 1$ is trivial and $\Sigma_n = S_n$ is the symmetric group, $\chi_{\textbf{1}, n}\defeq \textup{std}(S_n) $ is the standard $(n-1)$-dimensional irreducible character of $S_n$, inflated to $G\wr S_n$ via the epimorphism $G\wr S_n\twoheadrightarrow S_n$; explicitly $\forall v\in G^n, \sigma\in S_n: \quad \chi_{\textbf{1}, n}(v, \sigma) = \#\{\textup{fixed points of }\sigma\} - 1$. Thus when $\Sigma_n = S_n$, $\chi_{\phi, n}$ is a family of irreducible characters of dimension $\phi(1)n - \bbone_{\phi=1}$.

To state the main approximation result of the current paper about $G\wr S_n$, we need one more definition:

\begin{definition}
    Let $G$ be a compact group. We define $\Q_G\subseteq \C$ as a field extension of $\Q$, by adding all possible word measures on $G$ to $\Q$:
    \[ \Q_G\defeq \Q\left( \left\{\EX_w\left[ \phi \right] \right\}_{{\phi\in \hat{G},\,\, w\in F_r}}\right). \]
\end{definition}

When $G$ is finite, $\Q_G$ is contained in the minimal splitting field of $G$ (that is, the minimal field containing all complex character values of $G$), which is contained in the cyclotomic field extension $\Q\left[\textup{exp}\left({2\pi i} / {|G|}\right)\right]$. 
However, $\Q_G$ can be strictly smaller than the minimal splitting field: for instance, $\Q_{C_m} = \Q$ for every $m\in \N$.

For every word $w\in F_r$ and every basis $B\subseteq F_r$, we denote by $|w|_B$ the length of $w$ when written in basis $B$.
We also denote $|w| \defeq \min_B |w|_B$, where $B$ runs over all bases of $F_r$.

\begin{theorem}
\label{thm_main_result} (Main Result)
Let $w\in F_r$, let $G$ a compact group and let $\phi\in \hat{G}$ be an irreducible character.
Then for every $n\ge |w|$, the expectation $\EX_w\left[\chi_{\phi, n}\right]$ of $\chi_{\phi, n}\in \widehat{G\wr S_n}$ coincides with some rational function in $\Q_G(n)$, and
    \[ \EX_w\left[\chi_{\phi, n}\right] = \mathscr{C}_{\phi}(w) \cdot n^{1-\pi_{\phi}(w)} + O\left(n^{-\pi_{\phi}(w)}\right).\]
Moreover, if $\EX_w\left[\chi_{\phi, n}\right]\neq 0$ then $\EX_w\left[\chi_{\phi, n}\right] = \Omega\left(n^{-|w|}\right) $.
\end{theorem}
\noindent Note that this theorem generalizes \cite[Theorem 1.8]{PP15} and \cite[Theorem 1.11]{magee2021surface}.

\noindent In \cite[Conjecture 1.13]{HP22}, Hanany and Puder conjectured that for some families $\chi_n$ of irreducible characters of groups $G_n$, we have:
\begin{conjecture}
\label{conj_HP22}
For every word $w\in F_r$,
\[ \EX_w[\chi_n] = O\left(\dim\left(\chi_n\right)^{1-\pi(w)}\right).  \]
\end{conjecture}
\noindent The theorem above proves the conjecture in the case of $\chi_{\phi, n}.$

The families $(\chi_n)_{n\in\N}$ of characters mentioned in Conjecture~\ref{conj_HP22} are called stable characters.
In the subsequent paper \cite{Sho23II} we analyze word measures on stable characters of wreath products.
Informally, let $G$ be a finite group, and consider the ascending sequence of groups $G\wr S_{\bullet}\defeq (G\wr S_1\le G\wr S_2\le \ldots)$.
An irreducible stable character of $G\wr S_{\bullet}$ is a sequence $(\chi_n)_{n\in\N}$ of irreducible characters $\chi_n\in\widehat{G\wr S_n}$ that satisfies certain stability conditions, which we do not elaborate on here.
In \cite{Sho23II} we prove that for every word $w\in F_r$, for every finite group $G$, and for every irreducible stable character $(\chi_n)_{n\in \N}$ of dimension $\Omega(n^2)$ (equivalently, a sequence which is non-trivial and not one of the sequences $\chi_{\phi, n}$ considered in the current paper), we have
\begin{equation}
    \label{eq_main_res_of_ShoII}
    \EX_w[\chi_n] = O\left(n^{-\pi(w)}\right). 
\end{equation}
This generalizes the result of \cite{HP22} and represents a significant step towards Conjecture~\ref{conj_HP22}.

\subsection*{Additional Results}

\subsubsection*{Iterated Wreath Products}

This work shows that given a character $\phi$ of a compact group $G$ and a number $n$, we can get a new character $\textup{Ind}_n\phi$ of the group $G\wr S_n$, and analyze its $w$-measure in terms of $w$-measures of $\phi$. We also show that if $\phi\in \hat{G}-\{\textbf{1}\}$ is irreducible and non-trivial, then so is $\textup{Ind}_n\phi\in \widehat{G\wr S_n}-\{\textbf{1}\}$.
We can iterate this procedure: given $n_1, n_2$, $\textup{Ind}_{n_1, n_2}(\phi)\defeq \textup{Ind}_{n_2} \textup{Ind}_{n_1}\phi$ is a character of $(G\wr S_{n_1})\wr S_{n_2}$. In general we have:

\begin{definition}
\label{def_iterated_W_Ind}
($\mathcal{W}_{n_1, \ldots, n_m}(G), \textup{Ind}_{n_1, \ldots, n_m}\phi $)
Let $n_1, \ldots, n_m\in \N$, and let $G$ be a compact group with a character $\phi\in \textup{char}(G)$. 
Denote the iterated wreath product by
\begin{equation*}
    \begin{split}
        \mathcal{W}_{n_1, \ldots, n_m}(G) 
        &\defeq G\wr S_{n_1} \wr S_{n_2} \wr S_{n_3} \ldots \wr S_{n_m} \\
        &\defeq \left(\ldots\left(\left(G\wr S_{n_1}\right)\wr S_{n_2}\right)\wr S_{n_3}\ldots \right)\wr S_{n_m},
    \end{split}
\end{equation*}
and define a character $\textup{Ind}_{n_1, \ldots, n_m}\phi\in \textup{char}\left(\mathcal{W}_{n_1, \ldots, n_m}(G)\right)$ by induction on $m$, via 
\[\textup{Ind}_{n_1, \ldots, n_m} \phi\defeq \textup{Ind}_{n_m} \textup{Ind}_{n_1, \ldots, n_{m-1}} \phi.\]
Note that $\dim\left(\textup{Ind}_{n_1, \ldots, n_m} \phi\right) = \phi(1)\cdot n_1\cdots n_m$. 
We show (Proposition~\ref{prop_Ind_phi_is_irred}) that if $\phi\in \hat{G}-\{\textbf{1}\}$ then $\textup{Ind}_{n_1, \ldots, n_m} \phi\in \widehat{\mathcal{W}_{n_1, \ldots, n_m}(G)}-\{\textbf{1}\}$.
\end{definition}

It seems plausible that Conjecture~\ref{conj_HP22} may be generalized to many variables $n_1, \ldots, n_m$ instead of just one variable $n$, in view of the following theorem:

\begin{theorem}
\label{thm_main_iterated_result}
Let $G$ be a compact group, with a non-trivial irreducible character $\phi\in \hat{G}-\{\textbf{1}\}$, let $w\in F_r$ be a word and fix $m\in \N$.
For every sequence $n_1, \ldots, n_m\in \N$ satisfying $n_i \geq |w|$ for $i=1,\ldots, m$,
denote $\chi\defeq \textup{Ind}_{n_1, \ldots, n_m}\phi \in \widehat{\mathcal{W}_{n_1, \ldots, n_m}(G)}$.
Then $\EX_w\left[\chi\right]$ coincides with some rational function in $\Q_G(n_1, \ldots, n_m)$, and 
\[ \EX_w\left[\chi\right] = O\left(\dim(\chi)^{1-\pi(w)}\right) \]
where the implied constant depends on $w, m$ and $G$.
Moreover, if $\phi$ is a linear character, then 
\[ \EX_w\left[\chi\right] = (n_1\cdots n_m)^{1-\pi_{\phi}(w)}\cdot \left(C + O\left(\sum_{i=1}^m n_i^{-1}\right)\right), \]
for some positive integer $C\le |\textup{Crit}_{\phi}(w)|^m$ which depends on $w, m$.
\end{theorem}
\noindent For example, for $w = [x, y]$ we have $\EX_w\left[\chi\right] = \chi(1)^{-1} = \frac{1}{\dim(\phi)\cdot n_1\cdots n_m}$ (and also for $w = x^2 y^2$ if $\phi$ is a character of a real representation).

We also prove in Subsection~\ref{subsection_spherical_tree} the following corollary, which is quite surprising:

\begin{corollary}
\label{corollary_intro_spherical_tree}
Let $n_1, \ldots, n_m\in \N$, and let $\mathcal{T}$ be the following spherically symmetric rooted tree: there is a root with $n_m$ children, each child has $n_{m-1}$ children and so on (so there are $n_1\cdots n_m$ leaves). The automorphism group $\textup{Aut}(\mathcal{T}) = \mathcal{W}_{n_1, \ldots, n_m}(\{1\})$ acts on the set $\mathcal{T}_0$ of its leaves. 
It turns out that the permutation character $\#\textup{fix}\left(\textup{Aut}(\mathcal{T})\acts \mathcal{T}_0\right)$ of this action and the natural character $\#\textup{fix}\left(S_{n_m}\acts [n_m]\right)$ are very "close": for every word $w\in F_r$ of length at most $\min\{n_i\}_{i=1}^m$,
\begin{equation}
\label{equation_close_characters_of_tree}
    \left| \EX_w\left[\#\textup{fix}\left(\textup{Aut}(\mathcal{T})\acts \mathcal{T}_0\right)\right] - \EX_w\left[\#\textup{fix}\left(S_{n_m}\acts [n_m]\right)\right] \right| = O\left((n_{m-1} n_m)^{1-\pi(w)}\right).
\end{equation}
where the implied constant depends on $m, w$.
\end{corollary}

For example, by taking $n_{m-1} = n_m = n\to \infty$, the expectations in equation~\eqref{equation_close_characters_of_tree} are $1 + \Theta\left(n^{1-\pi(w)}\right)$ and the error term is $O\left(n^{2(1-\pi(w))}\right)$, which is much smaller. One can deduce, similarly to \cite{Expansion_Puder_2015} and \cite[Chapter 8]{HP22}, that the random Schreier graph obtained by connecting $\mathcal{T}_0$ by $r$ random symmetries of $\mathcal{T}$, thus making it a $2r$-regular graph, is asymptotically almost surely a close-to-optimal expander.

\subsubsection*{More About Primitivity Ranks}

The next result is about the $\textup{Aut}(F_r)$-invariants defined in Definition~\ref{def_witnesses_no_alg}: the $\phi$-primitivity ranks $\pi_{\phi}$.
\begin{theorem}
\label{thm_pi_phi_inequalities}
We combine some results about $\pi_{\phi}$:
\begin{enumerate}
    \item Proposition~\ref{prop_pi_phi_of_independent_words}: 
    Let $G$ be a compact group and let $\phi\in \hat{G}-\{\textbf{1}\}$.
    Let $w_1, w_2\in F_r$ be "disjoint" words: that is, there is a free decomposition $F_r = J_1 * J_2$ such that $w_i$ can be conjugated into $J_i$, for $i=1, 2$.
    Then
    \begin{align*}
         \pi_{\phi}(w_1w_2) & = \pi_{\phi}(w_1) + \pi_{\phi}(w_2), \\
         \mathscr{C}_{\phi}(w_1 w_2) & = \frac{\mathscr{C}_{\phi}(w_1) \mathscr{C}_{\phi}(w_2)}{\dim\phi}, \\
         \textup{Crit}_{\phi}(w_1 w_2) & = \{H_1 * H_2: \,\,H_i\in \textup{Crit}_{\phi}(w_i)\}. 
    \end{align*}

    This is a generalization of \cite[Lemma 6.8]{Puder_2014}.
    \item Corollary~\ref{corollary_pi_limit_std}: Denote by $\textup{std}\left(S_n\right)$ the standard irreducible character of $S_n$, that is $\textup{std} = \textup{fix} - 1$. Then for every $w\in F_r$,
    \[\pi(w) = \lim_{n\to \infty} \pi_{\textup{std}\left(S_n\right)}(w) = \min_{n\in \N} \pi_{\textup{std}\left(S_n\right)}(w). \]
    This corollary follows from \cite[Theorem (1.4’)]{PP15}.
    \item Corollary~\ref{corollary_pi_phi_of_nilpotent}: If $G$ is a finite $p$-group (for some prime $p$) and $\phi\in \hat{G}-\{\textbf{1}\} $, then for every word $w\in F_r$,
    \[\pi_{\phi}(w)\ge \pi_{C_p}(w).\]
    This corollary follows from \cite[Lemma 7]{COCKE2019440}.
\end{enumerate}
\end{theorem}

We have further results, which we mention here only briefly:
\begin{itemize}
    \item [App~\ref{section_general_wreath}:] 
    Recall Theorem~\ref{thm_main_result}: for every compact $G$, $\phi\in \hat{G}-\{\textbf{1}\}$ and $w\in F_r$, we consider the $w$-measure on $G\wr_{[n]} S_n$ and show that $\EX_w\left[\textup{Ind}_n\phi\right] = O\left(n^{1-\pi_{\phi}(w)}\right)$. In particular, if $w$ is not a proper power, then $\EX_w\left[\textup{Ind}_n\phi\right] = O\left(\dim(\textup{Ind}_n\phi)^{-1}\right)$.
    It is natural to ask - could we replace $S_n\acts[n]$ by some different sequence of group actions?
    
    Note that in the main result~\eqref{eq_main_res_of_ShoII} of \cite{Sho23II} and in Conjecture~\ref{conj_HP22}, $\textup{Ind}_n\phi$ is generalized to stable characters. This hints that a sequence of group actions that replaces $S_n\acts[n]$ should have some stable-like properties.
    
    In Theorem~\ref{thm_general_wreath_character_bound} we prove that every sequence of \textbf{asymptotically oligomorphic} group actions $\Sigma_n\acts X_n$ (see Definition~\ref{def_asymp_oligom}) may replace $S_n\acts[n]$: if $w$ is not a proper power, then $\EX_w\left[\textup{Ind}_{X_n}\phi\right] = O\left(\dim(\textup{Ind}_{X_n}\phi)^{-1/2}\right)$. We conjecture in Conjecture~\ref{conj_reiter_for_word_measure_on_ao_grps} that the $\frac{1}{2}$ can be replaced by $1$.

    The "asymptotically oligomorphic" condition is very simple.
    It generalizes the stability condition, and is applicable for every exhaustion of an oligomorphic group by finite groups.

    \item [App~\ref{appendix_free_product_of_cyclic}:] We prove a "wreath product" analogue of \cite[Conjecture 7.1]{PZimhony}:
    Let $\mathcal{F} = C_m * \ldots * C_m$ be the free product of $r$ copies of $C_m$.
    For every $\gamma\in \mathcal{F}$ and a compact group $G$ with $\phi\in \hat{G}-\{\textbf{1}\}$, we prove the bound
    \[\EX_{\alpha\sim U(\textup{Hom}(\mathcal{F}, G\wr S_n))} \left[\textup{Ind}_n \phi(\alpha(\gamma))\right] = O\left(n^{\chi^{\textup{max}}(\gamma)}\right)\]
    where $\chi^{\textup{max}}(\gamma)$ is defined in \cite[Conjecture 7.1]{PZimhony}. 
\end{itemize}

A fundamental tool in this paper is the "induction-convolution lemma", which plays a key role in almost all of the theorems above (see Section~\ref{section_ICL}).

\subsection*{Motivation}

\textbf{Invariants of Free Words: } \label{invariants_of_free_words}
    There is a plethora of quantitative properties of words in free groups. 
    Words have length (with respect to some specific basis), width (\cite{segal_2009}), height (\cite{Kim_2010}) and even thickness (\cite{Goldstein1999TheLA}). 
    We are interested in further properties. These properties, like some of the ones mentioned above, are invariant under the action of $\textup{Aut}(F_r)\acts F_r$, which make them relevant e.g.\ for identifying orbits of words (see the next item). 
    For example, instead of length in some basis, we can look at the shortest length in the automorphism orbit (and then a word has length 1 if and only if it is primitive).
    
    
    In Table~\ref{table:invariants_of_free_words}, we present some functions $F_r\to \C$ which are invariant under the action of $\textup{Aut}(F_r)\acts F_r$.
    For every such invariant function, we associate a sequence of compact groups that is related to the invariant: explicitly, some of its irreducible characters $\chi$ (the "stable" ones) tend (sometimes only conjecturally) to satisfy $\EX_w[\chi] = O\left(\dim(\chi)^{\textup{invariant}(w)}\right)$, and we refer to the relevant papers for explicit, exact formulations.
    \begin{table}[ht!]
    \centering
    \begin{tabular}{|| c | l | c | c | c ||} 
     \hline
     Invariant
     & Definition
     & Group
     & Papers \\ [0.5ex] 
     \hline\hline
     $\textup{cl}(w)$ & $\min\left\{g: \exists u_i, v_i: w = \prod_{i=1}^g [u_i, v_i] \right\}$ 
     & $U(n)$ 
     & \cite{unitary_Magee_2019}\\ [0.5ex] \hline
     $\textup{sql}(w)$ & $\min\left\{g: \exists u_i: w = \prod_{i=1}^g u_i^2 \right\}$ 
     & $O(n)$ & \cite{orthogonal_magee2019matrix}\\ [0.5ex] \hline
     $\textup{scl}(w)$ & $\lim_{n\to \infty} {\textup{cl}(w^n)}/{n}$ 
     & $U(n)$ & \cite{unitary_Magee_2019}\\ [0.5ex] \hline
     $\pi(w)$ & $\min\left\{\textup{rk}(H): H\le F_r, w\in H \textrm{ imprimitive}\right\}$ 
     & $S_n$ & \cite{PP15} \\ [0.5ex] \hline
     $\pi_{C_m}(w)$ & $\min\left\{\textup{rk}(H): H\le F_r, w\in K_m(H)\right\}$ 
     & $C_m\wr S_n$ & \cite{magee2021surface}\\ [0.5ex] \hline
     $\pi_q(w)$ & $\min\left\{\textup{rk}(\mathcal{I}): \mathcal{I}\triangleleft \F_q[F_r], 1-w\in \mathcal{I} \textrm{ imprimitive}\right\}$ 
     & $\textup{GL}_n(\F_q)$ & \cite{ernstwest2021word}\\ [0.5ex] \hline
     $\pi_{\phi}(w)$ & $\min\left\{\textup{rk}(H): \inner*{w} \lneqq_{\textup{alg}} H, \EX_{w\to H}[\phi] \neq 0\right\}$ 
     & $G\wr S_n$ & Current\\
     [1ex] 
     \hline
    \end{tabular}
    \caption{Invariants of Free Words. Here only, the notation $\triangleleft$ stands for a proper right ideal.}
    \label{table:invariants_of_free_words}
    \end{table}
    
    Theorem~\ref{thm_pi_phi_inequalities} joins some known inequalities between the invariants in Table~\ref{table:invariants_of_free_words}. 
    For every $w\in F_r$, 
    \begin{enumerate}
        \item $2\cdot \textup{scl}(w) \le 2\cdot \textup{cl}(w) - 1$.
        \item 
        $\pi(w) \le \pi_{\Ss^1}(w) \le 2\cdot \textup{cl}(w)$. 
        \item Based on numerical computer experiments, \cite[Conjecture 6.3.2]{nicolaus2019a} conjectured that $\pi(w) - 1\le 2\cdot \textup{scl}(w)$. This was conjectured also in \cite[Conjecture 1.14]{HP22}.
        \item 
        $\pi_q(w) \le \pi(w) $ (\cite[Proposition 1.8]{ernstwest2021word}), and it is conjectured that $\pi_q(w) = \pi(w)$ (\cite[Conjecture 1.9]{ernstwest2021word}).        
        \item For every non-trivial irreducible character $\phi$ of a compact group, $\pi_{\phi}(w)\ge \pi(w)$.
    \end{enumerate}
    
Let $p$ be a prime number and $m\in \N$ divisible by $p$.
By adding some non-trivial irreducible character $\phi_p$ of some $p$-group to the picture, we get an interesting graph of inequalities:

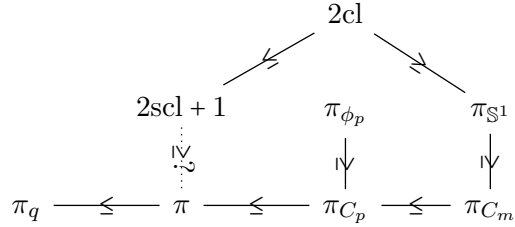
\[\begin{tikzcd}
	&& {2\textup{cl}} \\
	& {2\textup{scl}+1} & {\pi_{\phi_p}} & {\pi_{\mathbb{S}^1}} \\
	{\pi_q} & \pi & {\pi_{C_p}} & {\pi_{C_m}}
	\arrow["\ge"{marking}, no head, from=2-4, to=3-4]
	\arrow["\ge"{marking}, no head, from=2-3, to=3-3]
	\arrow["{\ge?}"{marking}, dotted, no head, from=2-2, to=3-2]
	\arrow["\le"{marking}, no head, from=3-1, to=3-2]
	\arrow["\le"{marking}, no head, from=3-3, to=3-4]
	\arrow["\le"{marking}, no head, from=2-2, to=1-3]
	\arrow["\ge"{marking}, no head, from=1-3, to=2-4]
	\arrow["\le"{marking}, no head, from=3-2, to=3-3]
\end{tikzcd}\]

\noindent \textbf{Distinguishing Orbits of Free Words:} \label{identify_free_word_up_to_auto} 
A conjecture by \cite[Question 2.2]{Amit2011CHARACTERSAS} and \cite[Conjecture 4.2]{Sha13}
asks whether every word in a free group is determined up to automorphism by its induced word measures on finite groups. 
(The converse, that two words in the same orbit induce the same measure on every finite group, is a simple observation). 
This conjecture has a slightly weaker version, replacing finite groups by compact groups.
This conjecture was proven for some specific (orbits of) words: 
\begin{itemize}
    \item Primitive Words (\cite{PP15}): $x$ (a basis element) is determined by the symmetric groups $S_n$.
    \item Surface words (\cite{magee2021surface}): $[x_1, y_1]\cdots [x_g, y_g]$ (orientable) and $x_1^2\cdots x_g^2$ (non-orientable) are determined by the unitary and orthogonal groups $U(n), O(n)$ and by wreath products $\Ss^1\wr S_n, C_2\wr S_n$.
    \item The commutator $[x, y]$ (\cite{Hanany_2020}).
    \item Powers of words that satisfy the strong (i.e.\ finite version) conjecture (\cite{Hanany_2020}): for example $x^d, [x, y]^d$.
\end{itemize}

In \cite{magee2021surface}, the invariants $\pi_{C_m}(w), \left|\textup{Crit}_m(w)\right|$ were used to recognize if $w$ is in the orbit of a surface word. 
Here we generalize these invariants to $\pi_{\phi}(w), \mathscr{C}_{\phi}(w)$ for every character $\phi$ of every compact group, which may lead to determination of more orbits of words.

\subsection*{Overview of the paper\footnote{Tip: for the readers convenience, this paper has a lot of links. When reading PDF with Adobe Reader, use Alt + $\leftarrow$ and Alt + $\rightarrow$ to jump backward and forward using the links.}}

\begin{enumerate}
    \item [Sec~\ref{section_morphisms_free_alg}:] We recall the definitions of free and algebraic extensions of free subgroups, and the algebraic-free decomposition (Theorem~\ref{thm_AFD}).
    We define the $B$-decompositions lattice $\textup{Decomp}_B(\eta)$ of a morphism $\eta$ (Definition~\ref{def_decomp_B}), with a focus on surjective morphisms $\textup{sur}_B, \mathcal{Q}_B(\Gamma)$ (Definition~\ref{def_surB_QB}), the algebraic version $\textup{Decomp}_{\textup{alg}}(\eta)$ (Definition~\ref{def_decomp_alg_B}) and the Möbius inversions on these lattices, especially the left Möbius inversion $L_{\eta}^B(n) $ (Proposition~\ref{prop_LB}).
    
    \item [Sec~\ref{section_ICL}:] We state the Induction-Convolution Lemma, which is the most fundamental tool in this paper, and some of its consequences. 
    We start with a basis-dependent formulation (Lemma~\ref{lemma_ICL}), and then deduce a basis-free formulation (Lemma~\ref{lemma_AlgICL}).
    
    \item [Sec~\ref{subsection_thm_main_result}:] We prove the main result of the paper, Theorem~\ref{thm_main_result}.
    \item [Sec~\ref{subsection_iterated_wreath}:] We apply the induction-convolution lemma recursively to prove Theorem~\ref{thm_main_iterated_result}. We deduce Corollary~\ref{corollary_intro_spherical_tree} about spherically symmetric trees.
    \item [Sec~\ref{section_witnesses}:] We investigate the invariants $\pi_{\phi}$, and prove Theorem~\ref{thm_pi_phi_inequalities}.
    We start with general groups, then continue to $p$-groups.
    
    \item [App~\ref{section_general_wreath}:] We combine the induction-convolution lemma with a theorem of Reiter (\cite{reiter19}) to deduce Theorem~\ref{thm_general_wreath_character_bound}, giving a bound on $\EX_w[\chi]$ in sequences of groups generalizing $G\wr S_n$.
    
    \item [App~\ref{appendix_free_product_of_cyclic}:] We prove a "wreath product analogue" of a conjecture from \cite{PZimhony}.
    \item [App~\ref{section_open_questions}:] We conclude with some open questions naturally arising from this work.
\end{enumerate}

\subsubsection*{Acknowledgments}
This paper is a part of an M.Sc.\ Thesis written by the author under the supervision of Prof.\ Doron
Puder from Tel Aviv University. This project has received funding from the European Research
Council (ERC) under the European Union’s Horizon 2020 research and innovation programme
(grant agreement No 850956).

I am deeply indebted to Professor Doron Puder, who guided me throughout the research, gave helpful advices for work and life in general, and taught me so much. This paper would not have been possible without him.
Many thanks to Asael Reiter for his curiosity and collaboration, to Yaron Brodsky for his smart advice, and to Alon Heller for a nice example.\\

\section{Core Graphs, Free vs.\ Algebraic, and Möbius Inversions}
\label{section_morphisms_free_alg}

We start by re-defining $\pi_{\phi}$ in a unified way.
There are situations where free groups $H \le J$ are "equivalent" with respect to word measures:

\begin{definition}
\label{def_free_factor}
(Free Factor) A subgroup $H\le F$ of a free group is called a \textbf{free factor} if there exists $J\le F$ such that $H * J = F$.
In this case we denote $H\leff{F}$, and say that $F$ is a \textbf{free extension} of $H$.
If $H\neq F$, $H$ ($F$) is called a \textbf{proper} free factor (extension, respectively). 
\end{definition}

It is easy to see that for every word $w\in H\leff J$ and for every compact group $G$, the word measures $\mu_{w\to H}, \mu_{w\to J}$ on $G$ are the same\footnote{Indeed, taking some basis of $J$ that contains a basis of $H$, every basis element is sent under a random homomorphism $J\to G$ to a random element, independent of the other basis elements. Thus the word measures $\mu_{w\to H}, \mu_{w\to J}$ are the same.}, so for every $\phi\in \textup{char}(G)$, $\EX_{w\to H}[\phi] = \EX_{w\to J}[\phi]$.

The following concept goes back to \cite{Takahasi1951NoteOC}, and was further studied in \cite{KAPOVICH2002608}, \cite{MVW07}:
\begin{definition}
\label{def_algebraic_extension}
(Algebraic Extension) 
Given $H\le J\le F_r$, we say that $J$ is an \textbf{algebraic extension} of $H$ if $H$ is not contained in any proper free factor of $J$, and then denote $H\le_{\textup{alg}} J$. If $H\neq J$, $J$ is called a \textbf{proper} algebraic extension. For a word $w\in H$ we say that $H$ is an algebraic extension of $w$ if $\inner*{w} \le_{\textup{alg}} H$.
Every f.g.\ $H$ has only a finite number of algebraic extensions (see propositions~\ref{prop_surj_extension_properties}, \ref{prop_alg_extension_properties}).
\end{definition}

Now we give the "correct" equivalent definitions for $\pi_{\phi}, \textup{Crit}_{\phi}$:

\begin{proposition}
\label{prop_equiv_def_pi_phi}
($\phi$-Witnesses of Imprimitivity)
Let $G$ be a compact group with an irreducible character $\phi\in \hat{G}$, and let $w\in F_r$ be a word in the free group of rank $r$.
We define the $\phi$-witnesses of $w$ as the set of subgroups 
    \[ \textup{Wit}_{\phi}(w) \defeq \{H\le F_r: \inner*{w} \lneqq H \textup{ is a proper algebraic extension, and } \,\EX_{w\to H}[\phi] \neq 0\}.\]
Then the $\phi$-primitivity rank of $w$ is the minimal rank in $\textup{Wit}_{\phi}(w)$:
    \[ \pi_{\phi}(w) = \min\{\textup{rk}(H): H\in \textup{Wit}_{\phi}(w)\}, \]
with the convention $min(\emptyset)=\infty$. 
The set of subgroups achieving the minimum is 
\[ \textup{Crit}_{\phi}(w) = \{H\in \textup{Wit}_{\phi}(w): \textup{rk}(H) = \pi_{\phi}(w)\}.\]
Moreover, $\pi(w)\le \pi_{\phi}(w)$ and if $\pi(w) = \pi_{\phi}(w)$ then $\textup{Crit}_{\phi}(w)\subseteq \textup{Crit}(w)$.
\end{proposition}

\begin{proof}
We need to show that taking subgroups of minimal rank among a large set of subgroups (all subgroups $H$ satisfying $\EX_{w\to H}[\phi]\neq 0$, or subgroups $H$ that contain $w$ but not as a primitive element in the case $\phi=\textbf{1}$) is the same as taking subgroups of minimal rank from the restricted set $\textup{Wit}_{\phi}(w)$, so it suffices to prove $\textup{Crit}_{\phi}(w)\subseteq \textup{Wit}_{\phi}(w)$.

Let $J\in \textup{Crit}_{\phi}(w)$. 
If $J$ is not an algebraic extension of $w$, there is a proper free factor $w\in H\leff J$. 
In particular $\textup{rk}(H) < \textup{rk}(J) = \pi_{\phi}(w)$.
In the case $\phi = \textbf{1}$, $w$ is not primitive in $H$ (since a free factor of a free factor is again a free factor), contradicting the minimality of $\pi(w)$.

In the case $\phi\neq \textbf{1}$, $\EX_{w\to H}[\phi] = \EX_{w\to J}[\phi] \neq 0$, contradicting the minimality of $\pi_{\phi}(w)$.
It remains to prove that $J$ is a proper extension, i.e.\ $J\neq \inner{w}$.
In the case $\phi = \textbf{1}$, it follows as $w$ is primitive in $\inner{w}$.
In the case $\phi\neq \textbf{1}$, $\mu_{w\to w}$ is the Haar (uniform) distribution on $G$, and by orthogonality of irreducible characters we have $\EX_{w\to \inner{w}} [\phi] = \inner*{\phi, \textbf{1}} = 0$.
Note that for every witness subgroup $H\in \textup{Wit}_{\phi}(w)$, $w$ is not primitive in $H$, as $H$ is a proper algebraic extension, hence $\pi(w)\le \pi_{\phi}(w)$.
\end{proof}

\subsection*{Core graphs}

This subsection is based on \cite{PP15} and \cite[Section 3.1]{HP22}. 
Let $B = \{b_1, \ldots, b_r\}$ be a basis of $F_r$, and consider the bouquet $\Omega_B$ on $r$ circles with distinct labels
from $B$ and arbitrary orientations and with wedge point $\otimes$. 
Then $\pi_1(\Omega_B, \otimes)$ is naturally identified with $F_r$.
For example,
\[\pi_1\left(\begin{tikzcd}
	\otimes \arrow[loop left]{l}{b_1} \arrow[loop right]{r}{b_2}
\end{tikzcd}, \,\otimes\right) = \textup{Free}(\{b_1, b_2\}).\]

The notion of ($B$-labeled) core graphs, introduced in \cite{Stallings1983TopologyOF}, refers to finite, connected, rooted graphs with every vertex (with a possible exception of the root) having degree at least two (so no leaves and no isolated vertices), that come with a graph morphism to $\Omega_B$ which is an immersion, namely, locally injective. 
In other words, this is a finite connected graph with at least one edge and no leaves (with a possible exception of the root), with edges that are directed and labeled by the elements of $B$, such that for every vertex $v$ and every $b\in B$, there is at most one incoming $b$-edge and at most one outgoing $b$-edge at $v$. 
By "rooted" we mean that there is a special vertex, the root, which is allowed to be a leaf.
We stress that multiple edges between two vertices and loops at vertices are allowed.

There is a natural one-to-one correspondence between finite $B$-labeled core graphs and f.g.\ non-trivial subgroups of $F_r$. 
Indeed, given a core graph $\Gamma$ as above, with a root $v$, 
consider the “labeled fundamental group” $\pi_1^{\textup{lab}}(\Gamma, v)$: closed paths in a graph with oriented and $B$-labeled edges correspond to words in the elements of $B$. 
In other words, if $p\colon \Gamma\to \Omega_B$ is the immersion, then $\pi_1^{\textup{lab}}(\Gamma, v)$ is the subgroup $p_*(\pi_1(\Gamma, v))$ of $\pi_1(\Omega_r, \otimes) = F_r$.

Conversely, if $H\le F_r$ is a f.g.\ subgroup, 
it corresponds to a finite core graph, denoted $\Gamma_B(H)$, which can be obtained in different manners. 
For example, let $\Upsilon$ be the rooted topological covering space of $\Omega_B$ corresponding to $H$, which is equal in this case to the Schreier graph depicting the action of $F_r$ on the right cosets of $H$ with respect to the generators $B$.
Then $\Gamma_B(H)$ is obtained from $\Upsilon$ by ’pruning all hanging trees’, or, equivalently, as the union of
all closed paths at the root, including the trivial constant path.
One can also construct $\Gamma_B(H)$ from any finite generating set of $H$ using “Stallings foldings” – see \cite{Stallings1983TopologyOF}, \cite{KAPOVICH2002608}, \cite{Puder_2014}, \cite{PP15} for more details about foldings and about core graphs in general.

\begin{example}
The core graph of the subgroup $H = \inner*{c, aca, a^{-1}ba}\le \textup{Free}(\{a, b, c\})$ is
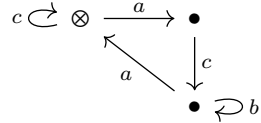
\[\begin{tikzcd}
	\otimes \arrow[loop left]{l}{c} & \bullet \\
	& \bullet \arrow[loop right]{r}{b}\\
	\arrow["a", from=1-1, to=1-2]
	\arrow["a", from=2-2, to=1-1]
	\arrow["c", from=1-2, to=2-2]
\end{tikzcd}\]
\end{example}

We cite from \cite{PP15} some basic facts (that were known before \cite{PP15}) about core graphs:

\begin{proposition}
\cite[Claim 3.1]{PP15} Let $H$ be a subgroup of $F_r$ with an associated core graph $\Gamma = \Gamma_B(H)$. 
The Euler Characteristic of a graph, denoted $\chi(\cdot)$, is the number of vertices minus the number
of edges. Then $rk(H) + \chi(\Gamma) = 1$. 
\end{proposition}

\begin{definition}
A \textbf{morphism} between core graphs is a map sending vertices to vertices, root to root, and edges to edges, that preserves the incidence relation and the orientations and labels of the edges.    
\end{definition}

\begin{proposition}
\label{prop_equivalence_of_categories}
\cite[Claim 3.2]{PP15}
Let $H, J, L\le F_r$ be subgroups. Then
\begin{enumerate}
    \item A morphism $\Gamma_B(H)\to \Gamma_B(J)$ exists if and only if $H\le J$, and if it
     exists, it is unique. We denote it by $\eta^B_{H\to J}$.
     In particular, whenever $H\le L\le J$, $\eta^B_{H\to J} = \eta^B_{L\to J}\circ \eta^B_{H\to L}$.
     \item Every morphism is an immersion (locally injective at the vertices).
     If $\eta^B_{H\to J}$ is injective, then $H\leff J$. 
\end{enumerate}
\end{proposition}

From now on, we will always have an ambient group $F_r$ of rank $r$ with a basis $B$ in the background.
The category of subgroups of non-trivial f.g.\ subgroups of $F_r$ with inclusions as morphisms is equivalent (by Proposition~\ref{prop_equivalence_of_categories}), via the functor $\Gamma_B$ (for any basis $B$), to the category of core graphs.
However, sometimes we want to be more specific, and make distinction between the words $w, w^{-1}$, even though they generate the same group:

\begin{definition}
Let $w\in H\le F_r$.
We define the morphism $\eta_{w\to H}$ as the data of the containment $w\in H$.
In contrast with $\eta_{\inner{w}\to H}$, a morphism $\eta_{w\to H}$ must "remember" the extra datum of the orientation of $w$, making the morphisms $\eta_{w\to H}, \eta_{w^{-1}\to H}$ different.
\end{definition}

 In graph-theoretic terms, whenever we have a core graph of a cyclic group $\inner*{w}$, we remember its orientation, and when having a morphism to another graph $\Gamma$ we remember the direction of the corresponding path.
This does not change a lot, and we will notify when this subtlety is required; for example, it is required in the following proposition, in which we give a geometric interpretation of $\EX_{w\to H}[\phi]$. For a core graph $\Delta$, denote by $E(\Delta)$ the set of topological edges, and by $\vec{E}(\Delta)$ the set of edges with orientations, so that $\left|\vec{E}(\Delta)\right| = 2\cdot \left| E(\Delta) \right|$.

\begin{proposition}
\label{prop_geometric_interpretation_of_E_eta}
Let $w\in J\le F_r$ be a word inside a subgroup $J$ of $F_r$ with core graph $\Delta\defeq \Gamma_B(J)$. 
Let $G$ be a compact group, and consider the set of anti-symmetric functions $\beta\colon \Vec{E}(\Delta)\to G$, that is $\beta(u\to v) = \beta(v\to u)^{-1}$; 
we can identify this space with $G^{E(\Delta)}$, thus we can draw a (Haar) uniformly random element $\beta$.
Then for every character $\phi\in \textup{char}(G)$,
\[ \EX_{w\to J}[\phi] = \EX_{\beta} \left[\prod_{e\in E(\Gamma_B(w))} \phi(\beta(\eta^B_{w\to J}(e)))\right], \]
where $\eta^B_{w\to J}(e)\in E(\Delta)$ is the image of $e$ under the graph-morphism $\eta^B_{w\to J}$, and the product is over the edges in the path $w$ in the order they appear in $w$.
\end{proposition}

\begin{proof}
We can fix a basis of $J$ by choosing a spanning tree $T\subseteq E(\Delta)$ and orientation for the rest of the edges:
for every $e\not\in T$ there is a basis element which is a path with end-points at the root and no back-tracking, that crosses $e$ once in the prescribed direction and avoids other edges which are not in $T$.
Regardless of the values $\{\beta(e)\}_{e\in T}$, the random variables $\{\beta(e)\}_{e\not\in T}$ are independent and Haar (uniform) distributed in $G$, thus defining a random homomorphism $\alpha \colon J\to G$. 
We finish the proof by observing that the $w$-path in $\Delta$ is evaluated as $\alpha(w)$ when we write $w$ in the basis corresponding to $T$.
\end{proof}

\subsection*{B-Surjective, Free and Algebraic Extensions}
\label{subsection_free_alg_ext}

Recall Definition~\ref{def_free_factor} of a free extension and Definition~\ref{def_algebraic_extension} of an algebraic extension.
Given a morphism $\eta\colon \Delta\to \Gamma$ of core graphs, we say that it is free or algebraic if $\pi_1^{\textup{lab}}(\Delta)$ is a free or algebraic extension of $\pi_1^{\textup{lab}}(\Gamma)$, respectively.

We cite the relevant definitions and theorems from \cite{PP15}.
For subgroups $H \le J\le F_r$, we denote by $\eta_{H\to J}\colon H\hookrightarrow J$ the inclusion map.
The following definition may seem unnecessarily complicated, as $\eta_{H\to J}$ is unique, but this formulation is useful for later.

\begin{definition}
\label{def_surB_QB}
($B$-surjective morphism, $\textup{sur}_B, \mathcal{Q}_B(\Gamma)$)
\cite[Definition 3.3]{PP15} Let $H\le J\le F_r$. Whenever $\eta^B_{H\to J}$ is surjective (on both vertices and edges), we say that $\Gamma_B(J)$ is a quotient graph of $\Gamma_B(H)$. We denote $H\Bcovers{J}$. The inclusion map $\eta_{H\to J}$ is then called $B$-surjective.
Define $\textup{sur}_B$ to be the collection of all $B$-surjective morphisms in the category of $B$-labeled core graphs.
Moreover, given a core graph $\Gamma$, define $\mathcal{Q}_B(\Gamma)$ as the set of morphisms in $\textup{sur}_B$ from $\Gamma$ to some quotient. 
\end{definition}

\begin{proposition}
\label{prop_surj_extension_properties}
\cite[Claim 3.4]{PP15}
Let $H, J, L\le F_r$ be f.g.\ subgroups. Then
\begin{itemize}
    \item If $H\Bcovers J$ and $H\le L\le J$, then $L\Bcovers J$.
    \item $H$ has only a finite number of $B$-quotient extension groups. 
    In particular, the poset of f.g.\ subgroups of $F_r$ with the partial order $\Bcovers$ is locally finite (i.e.\ $\forall x, z:\, |\{y: x < y < z\}| < \infty$).
\end{itemize}
\end{proposition}

Now we elaborate about algebraic morphisms.

\begin{proposition}
\label{prop_alg_extension_properties}
The following are some properties of algebraic morphisms: 
\begin{itemize}
    \item \cite[Claim 4.1]{PP15} The relation $\lealg$ is transitive. In other words, composition of algebraic morphisms is algebraic.
    \item \cite[Claim 4.2]{PP15} If $H\lealg J$ then for every basis $B\subseteq F_r$, the morphism $\Gamma_B(H)\to \Gamma_B(J)$ is surjective.
\end{itemize}
\end{proposition}

\begin{definition}
\label{def_alg(F_r)}
    Denote the set of all algebraic morphisms between subgroups of $F_r$ by $\textup{alg}(F_r)$. By Proposition~\ref{prop_alg_extension_properties}, for every basis $B\subseteq F_r$ we have $\textup{alg}(F_r)\subseteq \textup{sur}_B$.
\end{definition}

\begin{observation}
    \label{observe_alge_is_terminal_property}
    If $H_0 \le H_1 \le H_2$ and $H_0\to H_2$ is algebraic then $H_1\to H_2$ is algebraic too. 
    Indeed, if $H_1\le M \leff H_2$ then $H_0\le M\leff H_2$ implies $M=H_2$.
\end{observation}

\begin{theorem}
\label{thm_AFD}
(Free-Algebraic Decomposition; \cite[Claim 4.5]{PP15}) 
Let $H\le J$ be free groups. 
Then there is a unique subgroup $L$ of $J$ such that
$H\lealg L \leff J$.
Moreover, $L$ is the intersection of all intermediate free factors of $J$ and
the union of all intermediate algebraic extensions of $H$:
\[ L = \bigcap_{M: H\le M\leff J} M = \bigcup_{M: H\lealg J} M. \]
In particular, the intersection of all free factors is a free factor, and the union of
all algebraic extensions is an algebraic extension.
\end{theorem}

\subsection*{Möbius Inversions}
\label{subsection_mobius_inversions}

This subsection is based on \cite[Section 5]{PP15}. 
Throughout 
the paper we work with morphisms $\eta\colon H\to J$, which are just inclusion maps between subgroups of $F_r$.
While we could work with pairs of subgroups $(H,J)$ such that $H\leq J$ as in \cite{PP15}, we choose to use morphisms to ensure compatibility with a subsequent paper \cite{Sho23II}.

\begin{definition}
\label{def_decomp_B}
Let $\Gamma = \Gamma_B(H), \Delta = \Gamma_B(J)$ be core graphs of subgroups $H \le J$ such that $\eta\defeq \eta_{H\to J}$ is $B$-surjective.
Define
\begin{equation*}
    \begin{split}
        \textup{Decomp}_B(\eta) 
        &\defeq \{(\eta_1, \eta_2): \Gamma \overset{\eta_1}{\twoheadrightarrow} \Sigma \overset{\eta_2}{\twoheadrightarrow} \Delta \} \\
        &= \left\{(\eta_{H\to M}, \eta_{M\to J}): H \Bcovers M \Bcovers J \right\}.
    \end{split}
\end{equation*}
Similarly, let $\textup{Decomp}_B^3(\eta)$ denote the set of decompositions $\Gamma\overset{\eta_1}{\twoheadrightarrow}\Sigma_1\overset{\eta_2}{\twoheadrightarrow}\Sigma_2\overset{\eta_3}{\twoheadrightarrow}\Delta$ of $\eta$ into three surjective morphisms, and so on for $\textup{Decomp}_B^m(\eta)$ for all $m\in \N$. 
\end{definition}

Recall $\textup{sur}_B$ from Definition~\ref{def_surB_QB}.

\begin{definition}
($B$-Convolution)
Let $f, g \colon \textup{sur}_B\to \C$.
We define their $B$-convolution as follows: for every $B$-surjective morphism $\eta = \eta_{H \to J}$, 
\begin{equation*}
    \begin{split}
        (f \underset{B}{*} g)(\eta) 
        &\defeq \sum_{(\eta_1, \eta_2)\in \textup{Decomp}_B(\eta)} f(\eta_1) g(\eta_2)\\
        &= \sum_{H\Bcovers M \Bcovers J} f(\eta_{H\to M}) g(\eta_{M\to J}). 
    \end{split}
\end{equation*}
This operation, together with the obvious vector space structure, makes $\C^{\textup{sur}_B}$ into an associative convolution algebra.
The identity element is the function $\delta_B$ that assigns $1$ to isomorphisms and $0$ to any other morphism.
\end{definition}

\begin{definition}
(Möbius Inversion)
Let $1\in \C^{\textup{sur}_B}$ be the constant $1$ function. It is invertible\footnote{See for example the theory of incidence algebras in \cite{Sta97}. Alternatively, the restriction of the operator "convolution with 1" on a finite sub-poset acts as a uni-triangular matrix.} in $\C^{\textup{sur}_B}$, and we denote its inverse by $\mu^B$, called the Möbius inversion in the basis $B$. 
\end{definition}



Now we use Möbius inversions to define Möbius derivations. 
The following definition of $\EX\left[\Sigma_n\acts[n]\right]$ generalizes $\EX_{w\to H}[\phi]$, where $\phi$ is the permutation character of the group action $\Sigma_n\acts[n]$, to be defined over general $B$-surjective morphisms.
\begin{definition}
($\EX\left[\Sigma_n\acts[n]\right]$ and its Möbius Derivations)
Let $\Sigma_n\le S_n$ be a permutation group, acting on $[n]$.
We define
$\EX\left[\Sigma_n\acts[n]\right] \colon \textup{sur}_B\to \C$ as follows: given a morphism $\eta\colon H\to J$ (that is, $H\le J\le F_r$),
\[ \EX_{\eta}\left[\Sigma_n\acts[n]\right] \defeq \EX_{\alpha\sim U(\textup{Hom}(J, \Sigma_n))} \left[\#\textup{common fixed points of }\alpha(H)\le \Sigma_n \right]. \]
Also define its \textbf{L}eft, \textbf{R}ight and \textbf{C}entral Möbius derivations by
\[\begin{tikzcd}
	& {\mathbb{E}\left[\Sigma_n \curvearrowright[n]\right]} \\
	{\substack{L^B\left[\Sigma_n \curvearrowright[n]\right] \defeq \mu^B \underset{B}{*} \mathbb{E}\left[\Sigma_n \curvearrowright[n]\right]}} && {\substack{R^B\left[\Sigma_n \curvearrowright[n]\right] \defeq \mathbb{E}\left[\Sigma_n \curvearrowright[n]\right] \underset{B}{*} \mu^B}} \\
	& {\substack{C^B\left[\Sigma_n \curvearrowright[n]\right] \defeq \mu^B \underset{B}{*} \mathbb{E}\left[\Sigma_n \curvearrowright[n]\right] \underset{B}{*} \mu^B }}
	\arrow[from=1-2, to=2-3]
	\arrow[from=1-2, to=2-1]
	\arrow[from=2-1, to=3-2]
	\arrow[from=2-3, to=3-2]
\end{tikzcd}\]
that is, for every $H\Bcovers J$,
\begin{equation*}
\begin{split}
    \mathbb{E}_{H\to J}\left[\Sigma_n \curvearrowright[n]\right] \,
&= \sum_{H\Bcovers M \Bcovers J} L^B_{M\to J}\left[\Sigma_n \curvearrowright[n]\right]\, 
= \sum_{H\Bcovers M \Bcovers J} R^B_{H\to M}\left[\Sigma_n \curvearrowright[n]\right] \\
&= \sum_{H\Bcovers M_1 \Bcovers M_2 \Bcovers J} C^B_{M_1\to M_2}\left[\Sigma_n \curvearrowright[n]\right]. 
\end{split}
\end{equation*}
In the special case $\Sigma_n = S_n$, we adopt the notation from \cite[Section 5: Möbius inversions]{PP15} and denote 
$L^B(n)\defeq L^B\left[S_n \curvearrowright[n]\right], \,\,\,\,\,\,
R^B(n)\defeq R^B\left[S_n \curvearrowright[n]\right], \,\,\,\,\,\,
C^B(n)\defeq C^B\left[S_n \curvearrowright[n]\right].$
\end{definition}

We denote by $(n)_t \defeq n\cdot (n-1)\cdots (n-t+1)$ the falling factorial of $n$ of length $t$.

\begin{proposition}
("Basis dependent Möbius inversions",
\label{prop_LB}
(\cite[Lemma 6.4]{PP15})\\
Let $\eta\colon \Gamma\to \Delta$ be a $B$-surjective morphism of core graphs. Then for every $n\ge |E(\Gamma)|$,
\[ L^B_{\eta}(n) = \frac{\prod_{v\in V(\Delta)}(n)_{|\eta^{-1}(v)|}} {\prod_{e\in E(\Delta)}(n)_{|\eta^{-1}(e)|}}. \]
In particular,
\begin{enumerate}
    \item For every $n\ge |E(\Gamma)|$, $L^B_{\eta}(n)$ coincides with a rational function in $n$, with at most $|E(\Gamma)|$ poles (including multiplicity).
    \item $L^B_{\eta}(n) = n^{\chi(\Gamma)} \cdot (1 + O\left(n^{-1}\right))$.
    \item If there is only one letter (i.e.\ if $\Delta$ is a cycle) then $L^B_{\eta}(n)=1$.
\end{enumerate}
Note that $L^B_{\eta}[\Sigma_n\acts[n]]$ is generally much more involved.\footnote{This proposition can be slightly generalized: let $\Sigma_n\le S_n$. For every letter $b\in B$, denote by $E_b(\Gamma)$ the $b$-labeled edges of $\Gamma$, and $E_{max}(\Gamma)\defeq \max\{|E_b(\Gamma)|:b\in B\}$. 
If the action $\Sigma_n\acts[n]$ is $E_{\max}(\Gamma)$-transitive, then $L^B_{\eta}[\Sigma_n\acts[n]] = L^B_{\eta}(n)$.
Hence, if $\Sigma_n $ is $k$-transitive, then for every $w\in F_r$ with $E_{max}(w)\le k$, $\EX_w[\Sigma_n\acts[n]]=\EX_w[S_n\acts [n]]$.
For example, let $q$ be a prime power. Then $\textup{PGL}_2(\F_q) \acts \PR^1(\F_q)$ acts 3-transitively on the projective line over $\F_q$ of size $n = q + 1$, thus for every word $w\in F_r$ with $E_{\max}(w)\le 3$, $\EX_w\left[\textup{PGL}_2(\F_q) \acts \PR^1(\F_q)\right]$ is approximated by \cite{PP15}.
However, there are no $k$-transitive group actions for $k\ge 6$ except for the natural actions of $S_n, A_n\acts [n]$, and the set $B(k, r)\defeq \{w\in F_r : E_{max}(w)\le k\}$ of relevant words is small: $|B(k, r)|\le (2kr)^{r(k+1)}$ (it can be exactly computed, e.g.\ $|B(1, r)| = \floor{\sqrt{e}2^r r!}$). }
\end{proposition}

Similarly to the $B$-surjective case, we define the algebraic Möbius inversion and algebraic derivations, based on \cite[Definition 6.13]{HP22}. Recall $\textup{alg}(F_r)$ from Definition~\ref{def_alg(F_r)}.

\begin{definition}
\label{def_decomp_alg_B}
($\textup{Decomp}_{\textup{alg}}(\eta), \textup{Decomp}_{\textup{alg}}^3(\eta)$, Algebraic Möbius Inversion) 
For an algebraic morphism $\eta\colon \Gamma\to \Delta$
denote by $\textup{Decomp}_{\textup{alg}}(\eta)$ and $\textup{Decomp}_{\textup{alg}}^m(\eta)$ the set of decompositions of $\eta$ into two ($m$, respectively) algebraic morphisms.
If $\eta$ is algebraic, then $\textup{Decomp}_{\textup{alg}}(\eta)\subseteq \textup{Decomp}_B(\eta)$.
For every $f, g \colon \textup{alg}(F_r)\to \C$, we define their algebraic convolution by mapping every algebraic morphism $\eta = \eta_{H \to J}$ to
\[ (f \underset{\textup{alg}}{*} g)(\eta) \defeq \sum_{(\eta_1, \eta_2)\in \textup{Decomp}_{\textup{alg}}(\eta)} f(\eta_1) g(\eta_2). \]
Similarly to $\C^{\textup{sur}_B}$, $\C^{\textup{alg}(F_r)}$ is an associative convolution algebra, with the identity element $\delta_{\textup{alg}}$ that assigns $1$ to algebraic isomorphisms and $0$ to any other algebraic morphism.
As in the $B$-surjective case, the function $1\in \C^{\textup{alg}(F_r)}$ is invertible, and we denote its inverse by $\mu^{\textup{alg}}$.
\end{definition}

For the purposes of this paper, it is sufficient to define the algebraic Möbius derivatives of the group action $S_n\acts [n]$ only, rather than considering general permutation groups:
\begin{definition}
The algebraic left, right and central Möbius derivations of $\EX[S_n\acts [n]]$ are defined similarly to the $B$-surjective derivations, and are denoted by
$L^{\textup{alg}}_{\eta}(n), R^{\textup{alg}}_{\eta}(n), C^{\textup{alg}}_{\eta}(n)$, respectively. 
For instance, $L^{\textup{alg}}\defeq \mu^{\textup{alg}} \convalg \EX[S_n\acts [n]]$, or equivalently for every $\eta\in \textup{alg}(F_r)$,
\[ \EX_{\eta}[S_n\acts [n]] = \sum_{(\eta_1, \eta_2)\in \textup{Decomp}_{\textup{alg}}(\eta)} L_{\eta_2}^{\textup{alg}}(n). \]
\end{definition}

\section{The Induction-Convolution Lemma}
\label{section_ICL}

In this section we will handle the problem of computing $w$-expectation of an induced character. 
We will work with the following convenient equivalent definition of induced representation:
\begin{remark}
\label{remark_explicit_general_induction}
Generally, given a representation $\Phi\in \textup{Hom}(K, \textup{GL}_d(\C))$ of a finite index subgroup $K\le G$, one can compute the induced representation $\textup{Ind}_K^G \Phi$ explicitly with block-matrices.
Given a transversal $t_1, \ldots, t_n$ of $G/K$ and $g\in G$, the matrix $\textup{Ind}^G_K \Phi(g)$ has $n\times n$ blocks of size $d\times d$, where the $(i, j)$-block is
\begin{equation*}
    \textup{Ind}^G_K \Phi(g)_{i, j} = \begin{cases}
    \Phi(t_i^{-1} g t_j) &\textrm{if }t_i^{-1}gt_j\in K,\\
    0& \textrm{otherwise.}
    \end{cases}
\end{equation*}
\end{remark}

Recall the notation $\textup{Ind}_n\phi$ from Definition~\ref{def_general_Ind_phi}.

\begin{proposition}
\label{prop_Ind_phi_is_irred}
Let $G$ be a compact group, let $\phi\in\hat{G}-\{\textbf{1}\}$,
and let $\Sigma_n \le S_n$ be a transitive permutation group. 
Then $\textup{Ind}_n\phi\in \widehat{G\wr_{[n]} \Sigma_n }$ is an irreducible character.
\end{proposition}

\begin{proof}
Denote by $\textup{stab}(n)\le \Sigma_n $ the stabilizer of the point $n$. 
We start by showing that $\textup{Ind}_n\phi$ is actually an induced character.
Let $\Phi\colon G\to \textup{GL}_{\phi(1)}(\C)$ be the representation that yields $\phi$, that is $\textup{tr}\circ \Phi = \phi$.
We can extend this representation to $\Tilde{\Phi}\colon G\wr_{[n]} \textup{stab}(n)\to \textup{GL}_{\phi(1)}(\C)$ as follows: 
\[\forall v\in G^n, \sigma\in \textup{stab}(n):\quad \Tilde{\Phi}(v, \sigma)\defeq \Phi(v(n)).\] 
Define $\rho_{\phi}\colon G\wr_{[n]} \Sigma_n \to \textup{GL}_{\phi(1)n}(\C)$ that assigns to every $(v, \sigma)\in G\wr_{[n]} \Sigma_n $ the $n\times n$-block matrix whose $i, j$-block is $\Phi(v(i)) \cdot \bbone_{\sigma(i) = j}$.
This is a representation, induced from the extended representation $\Tilde{\Phi}$.
Indeed, choose representatives $t_1, \ldots, t_n\in \Sigma_n / \textup{stab}(n)$ of $(G\wr \Sigma_n)/(G\wr \textup{stab}(n))\cong \Sigma_n/\textup{stab}(n)$, satisfying $t_i(n) = i$, and identify them as $G\wr_{[n]} \Sigma_n $-elements by the embedding $\Sigma_n \cong \{1\}\wr_{[n]} \Sigma_n$. 
Then for every $(v, \sigma)\in G\wr_{[n]} \Sigma_n $ and $i, j \in [n]$, the $(i, j)$-block in the $nd$-dimensional matrix $\textup{Ind}_n\phi(v, \sigma)$ is $\bbone_{\sigma(i) = j}\cdot \Tilde{\Phi}(t_{j}^{-1}\cdot (v, \sigma) \cdot t_i)$ (note that $n\overset{t_i}{\to} i \overset{\sigma}{\to} j \overset{t_j^{-1}}{\to} n$ and $\Tilde{\Phi}$ is defined on $G\wr_{[n]} \textup{stab}(n)$) and 
$t_{j}^{-1}\cdot (v, \sigma) \cdot t_i = \left(v(t_i(k))_{k=1}^n, t_{j}^{-1}\sigma t_i\right)$ so the $(i, j)$-block is $\Phi(v(i))$. 
Hence $\rho_{\phi}=\textup{Ind}_n\Phi$ coincides with Definition~\ref{def_general_Ind_phi}.

Now we show irreducibility. Denote $ \Tilde{\phi}\defeq \textup{tr}\circ \Tilde{\Phi}$.
By Frobenius reciprocity, 
\begin{equation*}
    \begin{split}
        \inner*{\textup{Ind}_n\phi, \textup{Ind}_n\phi}_{G\wr_{[n]} \Sigma_n }
        &= \inner*{\left(\textup{Ind}_n\phi\right)\restriction_{G\wr_{[n]} \textup{stab}(n)}, \Tilde{\phi}}_{G\wr_{[n]} \textup{stab}(n)} \\
        &= \int_{g_1, \ldots, g_n\in G} \int_{\sigma\in \textup{stab}(n)} \left(\sum_{1\le i\le n: \,\sigma(i)=i} \phi(g_i)\right) \overline{\phi(g_n)}\\
        &= \int_{ g_n\in G } \overline{\phi(g_n)} \int_{g_1, \ldots, g_{n-1}\in G} \int_{\sigma\in \textup{stab}(n)} \left(\phi(g_n) + \sum_{1\le i \le n-1: \,\sigma(i)=i} \phi(g_i)\right) \\
        &= \inner*{\phi, \phi}_G + \inner*{\textbf{1}, \phi}_G \cdot \inner*{\textup{Ind}_{n-1} \phi, \textbf{1}}_{G\wr_{[n-1]} \textup{stab}(n)}.
    \end{split}
\end{equation*}
Since $\phi\in\hat{G}-\{\textbf{1}\}$, we have $\inner*{\phi, \phi}_G = 1, \inner*{\textbf{1}, \phi}_G = 0$ so $\textup{Ind}_n\phi\in \widehat{G\wr_{[n]} \Sigma_n }$ is irreducible.
\end{proof}

\begin{remark}
In the case $\phi=\textbf{1}$, $\textup{Ind}_n\phi$ is just the permutation character of the action $\Sigma_n\acts[n]$, so in general $\chi_{\textbf{1}, n} = \left(\textup{Ind}_n\textbf{1}\right) - \textbf{1}$ is irreducible if and only if the action is doubly-transitive.
Indeed, the action $\Sigma_n\acts[n]$ is isomorphic to the action of $G\wr \Sigma_n $ on the $n$ cosets of $G\wr \textup{stab}(n)$.
Since $\EX_{\eta}[\Sigma_n\acts[n]]$ (which we defined for all morphisms $\eta\colon H\to J$, i.e.\ inclusions $H\le J$) generalizes $\EX_{\eta}[\textup{Ind}_n\textbf{1}]$ (which we defined only for $\eta\colon \inner{w}\to H$ for a word $w\in H$), we are led to denote $\EX_{\eta}[\textup{Ind}_n\textbf{1}] = \EX_{\eta}[\Sigma_n\acts[n]]$ for every $\eta$.
\end{remark}

Given $\phi\colon G\to \C$, we aim to break down the expression $\EX_{\eta}[\textup{Ind}_n\phi]$ into components that exclude the wreath product $G\wr_{[n]} \Sigma_n $, and instead only involve $G$ and the action $\Sigma_n\acts[n]$ separately.
The induction-convolution lemma gives such a decomposition, in the form of convolution in the lattice of $B$-surjective morphisms.

Recall $\mathcal{Q}_B(\Gamma)$ from Definition~\ref{def_surB_QB}; we define $\mathcal{Q}_B(w)$ as $\mathcal{Q}_B(\Gamma_B(\inner{w}))$ with additional data.
\begin{definition}
Let $w\in H\le F_r$.
For every basis $B\subseteq F_r$, we define $\mathcal{Q}_B(w)$ as the set of $B$-surjective quotient maps $\Gamma_B(w)\twoheadrightarrow \Gamma_B(H)$. 
By abuse of notation, we also write $H\in \mathcal{Q}_B(w)$ (for $H\in F_r$) or $\Gamma\in \mathcal{Q}_B(w)$ (for $\Gamma = \Gamma_B(H)$) instead of $\eta_{w\to H}\in \mathcal{Q}_B(w)$. 
\end{definition}

\begin{lemma}
\label{lemma_ICL}
(The Induction-Convolution Lemma) 
Let $B\subseteq F_r$ be a basis, and $w\in F_r$. 
Then
\[ \EX_w\left[\textup{Ind}_n\phi\right] = \sum_{H\in \mathcal{Q}_B(w)} \EX_{w\to H}[\phi]\cdot L_{H\to F_r}^B(\Sigma_n\acts[n]). \]
Another useful formulation of the lemma is the following equality between operators\footnote{Note that the convolution is well-defined even though $\EX[\phi], \EX[\textup{Ind}_n\phi]$ are defined only on $\eta\colon \inner{w}\to H$.
However, $\mu^B$ and $\EX[\Sigma_n\acts[n]] = \EX[\textup{Ind}_n\textbf{1}]$ are defined on every morphism, and this is crucial as convolution uses many morphisms that arise in $\textup{Decomp}_B(\cdot)$. }:
\[ \EX[\textup{Ind}_n\phi] = \EX[\phi] \underset{B}{*} \mu^B \underset{B}{*} \EX[\textup{Ind}_n\textbf{1}]\]
where $\EX[\phi], \EX[\textup{Ind}_n\phi]$ are elements of $\C^{\mathcal{Q}_B({w})}$: they map the morphism $\eta\colon w\to H$ to
$\EX_{\eta}[\phi], \EX_{\eta}[\textup{Ind}_n\phi]$ respectively. 
Equivalently, for every morphism $\eta\in \mathcal{Q}_B(w)$,
\begin{equation*}
\label{eq_ICL_withL}
    \EX_{\eta}\left[\textup{Ind}_n\phi\right] = \sum_{(\eta_1, \eta_2)\in \textup{Decomp}_B(\eta)} \EX_{\eta_1}[\phi]\cdot L_{\eta_2}^B(\Sigma_n\acts[n]). 
\end{equation*}
\end{lemma}

Before we prove the induction-convolution lemma, we describe some special cases that were proven before.
The simplest case of the induction-convolution lemma is when $\phi=\textbf{1}$.
In this case, the statement of the lemma is just the definition of the Möbius inversion $\mu^B$ as the inverse of the constant function $1$ in the convolution algebra $\C^{\textup{sur}_B}$. 
This is the case which is needed for \cite{PP15}, so one can think of this lemma as a way of lifting results from $S_n$ to $G\wr S_n$ when $G$ is not trivial.

For the next example, recall Definition~\ref{def_K_m(H)} ($K_m(H)$).
A more interesting case of the induction-convolution lemma was proved in \cite[subsection 3.1]{magee2021surface}. 
In this case, $\phi$ was an embedding of $\Z/m\Z$ into $\Ss^1$ (or the identity on $\Ss^1$ if $G=\Ss^1$), and 
$\EX_{\eta}[\phi] = \bbone_{w\in K_m(\textup{Im}(\eta))}$ (recall Example~\ref{equation_expectation_of_cyclic_embedding}): for $w\in J$,
\begin{equation*}
\EX_{w\to J}[\phi] =
    \begin{cases}
    1 & \textrm{if } w \in K_m(J),\\
    0 & \textrm{otherwise.}
    \end{cases}
\end{equation*}
Equivalently, $\EX_{w\to J}[\phi]$ is the indicator of the event that for each edge $e$ in the graph $\Gamma_B(J)$, the (signed) number of times that the $\eta_{*}(w)$-path crosses $e$ is $0$ mod $m$. (If $m=\infty$, it means the edge is crossed the same number of times in each direction).
If this happens we say that $\eta$ is $m$\textbf{-balanced}; then the induction-convolution lemma reads 
\begin{equation*}
    \EX_{\eta}[\textup{Ind}_n\phi] = \sum_{\substack{(\eta_1, \eta_2)\in \textup{Decomp}_B(\eta):\\ \eta_1 \textrm{ is }m\textrm{-balanced}}} L_{\eta_2}^B(n).
\end{equation*}

\begin{definition}
(Block-Trace)
Let $A\in M_{n\times n}(\F)$ be a matrix.
For every $d \divides n$ let
\[ tr^n_d(A) \defeq \sum_{i=1}^{n/d} A_{(i, i)} \in M_{d\times d}(\F) \]
where $A_{(i, i)}$ is the $(i, i)^{th}$ $d\times d$ block in $A$.
Clearly $tr = tr^n_1$ and if $a \divides b \divides n $ then $tr^b_a \circ tr^n_b = tr^n_a$.
\end{definition}

Now we are ready to prove the induction-convolution lemma:

\begin{proof}[Proof of Lemma~\ref{lemma_ICL}]
Let $w = b_{t(0)}^{\varepsilon(0)}\cdots b_{t(L-1)}^{\varepsilon(L-1)}$ for $\varepsilon\colon \Z/L \to \{\pm 1\}$ and $b_{t(i)}\in B = \{b_1, \ldots, b_r\}$, and let $\Phi\colon G\to \mathcal{U}(d)\le \textup{GL}_d(\C)$ be a unitary representation yielding $\phi$. 
Recall $\textup{Ind}_n\Phi\colon G\wr_{[n]} \Sigma_n \to \mathcal{U}(d\cdot n)$ from Definition~\ref{def_general_Ind_phi}.
Then
\begin{equation*}
    \begin{split}
     \EX_w\left[\textup{Ind}_n\phi\right] 
     &= \EX_{v_1,\ldots, v_r\sim U(G^n)} \EX_{\sigma_1, \ldots, \sigma_r\sim U(\Sigma_n) } \textup{tr}^{nd}_1\left( \textup{Ind}_n\Phi(w((v_1, \sigma_1), \ldots, (v_r, \sigma_r))) \right) \\
     &= \EX_{v_1,\ldots, v_r\sim U(G^n)} \EX_{\sigma_1, \ldots, \sigma_r\sim U(\Sigma_n) } \textup{tr}^d_1 \textup{tr}^{nd}_d\left( \prod_{\ell\in \Z/L} \textup{Ind}_n\Phi\left((v_{t(\ell)}, \sigma_{t(\ell)})^{\varepsilon(\ell)}\right) \right) \\
     &= \EX_{v_1,\ldots, v_r\sim U(G^n)} \EX_{\sigma_1, \ldots, \sigma_r\sim U(\Sigma_n) } \textup{tr}^{d}_{1} \left( \sum_{i\colon \Z/L \to [n]} \prod_{\ell\in \Z/L} \left[\textup{Ind}_n\Phi\left((v_{t(\ell)}, \sigma_{t(\ell)})^{\varepsilon(\ell)}\right)\right]_{i(\ell), i(\ell + 1)} \right)\\
     &= \sum_{i\colon \Z/L \to [n]} \EX_{v_1,\ldots, v_r\sim U(G^n)} \EX_{\sigma_1, \ldots, \sigma_r\sim U(\Sigma_n) } \textup{tr}^{d}_{1} \left(  \prod_{\ell\in \Z/L} \left[\textup{Ind}_n\Phi\left((v_{t(\ell)}, \sigma_{t(\ell)})^{\varepsilon(\ell)}\right)\right]_{i(\ell), i(\ell + 1)} \right)\\
     &= \sum_{i\colon \Z/L \to [n]} \EX_{v_1,\ldots, v_r\sim U(G^n)} \EX_{\sigma_1, \ldots, \sigma_r\sim U(\Sigma_n) } \textup{tr}^{d}_{1} \left(  \prod_{\ell\in \Z/L} \bbone_{\sigma_{t(\ell)}^{\varepsilon(\ell)}(i(\ell)) = i(\ell + 1)} \Phi\left(\Tilde{v}_{t(\ell)}^{\varepsilon(\ell)}(i(\ell))\right)\right)
    \end{split}
\end{equation*}
where the indices are block-indices, and $\Tilde{v}$ is either $v$ or $\sigma_{t(\ell)}^{-1}.v$ depending on $\varepsilon(\ell)$:
\begin{equation*}
    \Tilde{v}_{t(\ell)}^{\varepsilon(\ell)}(i(\ell)) = \begin{cases}
        v_{t(\ell)}(i(\ell))        & \textup{ if } \varepsilon(\ell) = 1,\\
        v_{t(\ell)}(i(\ell+1))^{-1} & \textup{ if } \varepsilon(\ell) = -1.\\
    \end{cases}
\end{equation*}
Now if we think of $\Z/L$ as the vertices of $\Gamma_B(w)$, then we can think of $i\colon V(\Gamma_B(w)) \to [n]$ as a coloring and then glue together all vertices of the same color to get a quotient labeled rooted graph $\Gamma$.
For example, if $w = [b_1, b_2]$ and $(i(0), \ldots, i(3)) = (3, 7, 7, 5)$ then the $w$-path and the quotient graph are

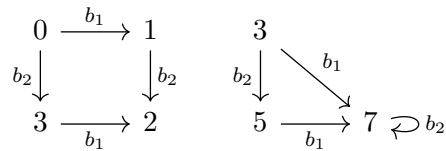
\[\begin{tikzcd}
	0 & 1 & 3 \\
	3 & 2 & 5 & 7 \arrow[loop right]{l}{b_2}
	\arrow["{b_1}", from=1-1, to=1-2]
	\arrow["{b_2}", from=1-2, to=2-2]
	\arrow["{b_2}"', from=1-1, to=2-1]
	\arrow["{b_1}"', from=2-1, to=2-2]
	\arrow["{b_2}"', from=1-3, to=2-3]
	\arrow["{b_1}", from=1-3, to=2-4]
	\arrow["{b_1}"', from=2-3, to=2-4]
\end{tikzcd}\]
Then we can split the sum according to this quotient graph, now summing only over injective functions
\[ \sum_{\Gamma} \sum_{i\colon V(\Gamma)\hookrightarrow [n]} \EX_{v_1,\ldots, v_r\sim U(G^n)} \EX_{\sigma_1, \ldots, \sigma_r\sim U(\Sigma_n) } \textup{tr}^{d}_{1} (\ldots) \]
and observe that the factor $\bbone_{\sigma_{t(\ell)}^{\varepsilon(\ell)}(i(\ell)) = i(\ell + 1)}$ guarantees that the contribution of a graph $\Gamma $ is 0 if it is not a core graph: if a vertex $v = i(\ell) \in [n]$ has, say, two $b_1$ outgoing edges, then both of the target vertices must be $\sigma_1(v)$, contradicting injectivity of $i$.
In the example above, the equations $\sigma_1(3) = \sigma_{t(0)}^{\varepsilon(0)}(i(0)) = i(1) = 7, \sigma_1^{-1}(7) = \sigma_{t(2)}^{\varepsilon(2)}(i(2)) = i(3) = 5$ cannot hold together so the contribution of this "non-core" graph vanishes.
So the sum runs over $\Gamma\in \mathcal{Q}_B(w)$. Let us call a function $i\colon V(\Gamma_B(w)) \hookrightarrow [n]$ \textbf{valid} if for every $b$-labeled edge $u_1\overset{b}{\to} u_2$, we have $\sigma_b(i(u_1)) = i(u_2)$. Then 
\begin{equation*}
    \begin{split}
     \EX_w\left[\textup{Ind}_n\phi\right] 
     &= \sum_{\Gamma\in \mathcal{Q}_B(w)} \sum_{i\colon V(\Gamma) \hookrightarrow [n]} \EX_{\substack{v_1,\ldots, v_r\sim U(G^n)\\\sigma_1, \ldots, \sigma_r\sim U(\Sigma_n) }} tr^{d}_{1} \left( \bbone_{i \textrm{ is valid}}  \prod_{e\in w\textrm{-path from root}} \Phi\left(v_{\textrm{label}(e)}(i(\textrm{src}(e)))\right)\right)
    \end{split}
\end{equation*}
where $\textrm{label}(e)\in \{b_1, \ldots, b_r\}$ is the label of $e$, 
$\textrm{src}(e)\in V(\Gamma)$ is the source vertex of $e$,
and the product runs over edges in the $w$-path in $\Gamma$ starting from the root of $\Gamma$, so there are exactly $L$ terms in the product, arranged according to their order in $w$.
Moreover, by substitution of $\phi = \textbf{1}$ (the trivial representation) we see that 
\begin{equation*}
    \begin{split}
     \EX_w\left[\textup{Ind}_n\textbf{1}\right] 
     &= \sum_{\Gamma\in \mathcal{Q}_B(w)} \sum_{i\colon V(\Gamma) \hookrightarrow [n]} \EX_{\substack{v_1,\ldots, v_r\sim U(G^n)\\\sigma_1, \ldots, \sigma_r\sim U(\Sigma_n) }}\bbone_{i \textrm{ is valid}} 
    \end{split}
\end{equation*}
which means that 
\[ L_{\pi_1(\Gamma)\to F_r}^B(\Sigma_n\acts[n]) = \sum_{i\colon V(\Gamma) \hookrightarrow [n]} \EX_{\sigma_1, \ldots, \sigma_r\sim U(\Sigma_n) }\bbone_{i \textrm{ is valid}}. \]
(We removed the expectation with respect to $v_1, \ldots, v_r\sim U(G^n)$ as they are not relevant for the validity of $i$).
Now we use the structure of the wreath product for the first time: 
All the $nr$ random variables $\left(v_{j}(k)\right)_{1\le j\le r, 1\le k\le n}$ are independent and Haar (uniform) distributed in $G$.
Since $i$ is injective, the random variables $\left\{v_{\textrm{label}(e)}(i(\textrm{src}(e)))\right\}_{e\in E(\Gamma)}$ are all distinct: if we have two edges $e_1, e_2$ with $\textrm{label}(e_1) = \textrm{label}(e_2), i(\textrm{src}(e_1)) = i(\textrm{src}(e_2)) $ then $\textrm{src}(e_1) = \textrm{src}(e_2) $ and as the graph $\Gamma$ is a core graph we must have $e_1 = e_2$.
This means that $v_{\textrm{label}(e)}(i(\textrm{src}(e)))$ are all independent and Haar (uniform) distributed in $G$, and by Proposition~\ref{prop_geometric_interpretation_of_E_eta}, we get
\begin{equation*}
\begin{split}
    \forall \Gamma\in \mathcal{Q}_B(w): \forall i\colon V(\Gamma) \hookrightarrow [n]:\quad 
    &\EX_{v_1,\ldots, v_r\sim U(G^n)} \left[ tr^{d}_{1} \prod_{e\in w\textrm{-path from root}} \Phi\left(v_{\textrm{label}(e)}(i(\textrm{src}(e)))\right)\right] \\
    =& \EX_{\beta\sim U(E(\Gamma)\to G)} \left[ tr^{d}_{1} \prod_{e\in w\textrm{-path from root}} \Phi\left(\beta(e)\right)\right] \\
    =& \EX_{\beta\sim U(E(\Gamma)\to G)} \left[ \phi\left( \prod_{e\in w\textrm{-path from root}} \beta(e)\right)\right] \\
    =& \EX_{w\to \pi_1(\Gamma)}[\phi]
\end{split}
\end{equation*}
 and this finishes the proof.
\end{proof}

\subsection{Algebraic Induction-Convolution Lemma}

We will soon state an algebraic version of the same lemma, but first we have to understand the connection between $L^{\textup{alg}}$ and $L^B$.

It follows from \cite[Proposition 6.15]{HP22} that for an algebraic morphism $\eta$,
\[ L^{\textup{alg}}_{\eta} = \sum_{\substack{(\eta_1, \eta_2)\in \textup{Decomp}_B(\eta):\\ \eta_1 \textrm{ is free}}} L^B_{\eta_2}. \]
\begin{corollary}
\label{corollary_L_alg_approx}
For every algebraic morphism $\eta\colon \Gamma\to \Omega_r$,
\[ L^{\textup{alg}}_{\eta} = n^{\chi(\Gamma)} \cdot \left(1 + O\left(n^{-1}\right)\right). \]
\end{corollary}

\begin{proof}
Choose some basis $B$; since for every $\eta'\colon \Gamma'\to \Omega_r$ we have $L^B_{\eta'} = n^{\chi(\Gamma')} \cdot \left(1 + O\left(n^{-1}\right)\right)$,
\begin{equation*}
    \begin{split}
    L^{\textup{alg}}_{\eta} (n)
    &= \sum_{\substack{(\eta_1, \eta_2)\in \textup{Decomp}_B(\eta):\\ \eta_1 \textrm{ is free}}} L^B_{\eta_2}(n)\\
    (\textup{Proposition}~\ref{prop_LB})\quad\quad\quad\quad &= \sum_{\substack{(\eta_1, \eta_2)\in \textup{Decomp}_B(\eta):\\ \eta_1 \textrm{ is free}}} n^{\chi(\textup{Im}(\eta_1))} \cdot \left(1 + O\left(n^{-1}\right)\right)\\
    &= n^{\chi(\Gamma)} \cdot \left(1 + O\left(n^{-1}\right)\right)
    \end{split}
\end{equation*}
where the last step is obtained by splitting the summation to the case where $\eta_1$ is an isomorphism, which contributes the dominant term, and every free morphism that is not an isomorphism, which must decrease Euler characteristic\footnote{Indeed, let $H\leff J$ be a free extension of subgroups with $\textup{rk}(J) \le \textup{rk}(H)$; then $H=J$.}.
\end{proof}

However, we prove a slightly more general version of this connection. We formulate everything in the language of morphisms, as this will be helpful for the subsequent paper \cite{Sho23II}.

\begin{definition}
\label{def_free-inv_func}
(Free-Invariant Function)
A function $g\in \C^{\textup{sur}_B}$ is called \textbf{free-invariant} if it is invariant with respect to post-composition with a free morphism, i.e.\ if 
$\eta_{0\to 2} = \eta_{1\to 2}\circ \eta_{0\to 1}$
and $\eta_{1\to 2}$ is free, as in the following diagram
\[\begin{tikzcd}
	{\Gamma_0} & {\Gamma_1} & {\Gamma_2}
	\arrow[from=1-1, to=1-2]
	\arrow["free", from=1-2, to=1-3]
\end{tikzcd}\]
then $g(\eta_{0\to 1}) = g(\eta_{0\to 2})$.\\
It is a standard fact that if a function depends only on the distribution of words under random homomorphisms to compact groups, then it is free-invariant. 
Indeed, if $F=H*J$ is a free decomposition, and $G$ is any group, then there is a natural identification
$\textup{Hom}(F, G) \cong \textup{Hom}(H, G) \times \textup{Hom}(J, G) $.
In particular if $G$ is compact and $\alpha\sim U(\textup{Hom}(F, G))$ is a random homomorphism, then $\alpha\restriction_H\sim U(\textup{Hom}(H, G))$ has Haar (uniform) distribution.
\end{definition}

Recall from Proposition~\ref{prop_alg_extension_properties} that algebraic morphisms are $B$-surjective. 

\begin{lemma}
\label{lemma_mobius_alg_vs_B}
($\mu^{\textup{alg}}$ vs. $\mu^B$)
Let $f\in \C^{\textup{sur}_B}$, and let $\eta$ be an algebraic morphism. 
Then
\begin{enumerate}
    \item 
    \[ (\mu^{\textup{alg}} \convalg f)(\eta) = \sum_{\substack{(\eta_1, \eta_2)\in \textup{Decomp}_{B}(\eta):\\ \eta_1 \textrm{ is free}}} (\mu^B \Bconv f)(\eta_2). \]
    In particular, since the operators $\Bconv, \convalg$ share the same identity element, we have
    \[ \mu^{\textup{alg}}(\eta) = \sum_{\substack{(\eta_1, \eta_2)\in \textup{Decomp}_{B}(\eta):\\ \eta_1 \textrm{ is free}}} \mu^B (\eta_2).\]
    \item For every free-invariant $g\in \C^{\textup{sur}_B}$, 
    \[ \left(g \Bconv \mu^B \Bconv f\right)(\eta) = \left(g \convalg \mu^{\textup{alg}} \convalg f\right)(\eta). \]
\end{enumerate}
\end{lemma}

\begin{proof}
Denote $\ell_B \defeq \mu^B\Bconv f $ and  $\ell_{\textup{alg}} \defeq \mu^{\textup{alg}}\convalg f $.
During the proof, the reader should have in mind the commutative diagram
\[\begin{tikzcd}
	& {\Gamma_1} \\
	{\Gamma_0} && {\Gamma_2} & {\Gamma_3}
	\arrow["alg"{description}, from=2-1, to=1-2]
	\arrow["free"{description}, from=1-2, to=2-3]
	\arrow["{\ell_B}"', from=2-3, to=2-4]
	\arrow["{\ell_{\textup{alg}}}", dashed, from=1-2, to=2-4]
	\arrow["g"', from=2-1, to=2-3]
\end{tikzcd}\]
with $\forall 0\le i < j \le 3: \,\eta_{i\to j}\colon \Gamma_i\to \Gamma_j$, where $\Gamma_0, \Gamma_3$ are fixed and $\Gamma_1, \Gamma_2$ vary, and regard the diagram (together with the hints above the arrows) as an informal illustrative description of the proof.
We will prove the first claim for $\eta = \eta_{1\to3}$, and prove the second claim for $\eta_{0\to 3}$.
For every $f\in \C^{\textup{sur}_B}$ and free-invariant $g\in \C^{\textup{sur}_B}$, 
\begin{equation}
\label{eq_proof_of_lemma_mobius_alg}
    \begin{split}
        (g\Bconv\mu^B\Bconv f)(\eta_{0\to 3})
        &= (g\Bconv \ell_B)(\eta_{0\to 3})\\
        &= \sum_{(\eta_{0\to 2}, \eta_{2\to 3})\in \textup{Decomp}_B(\eta_{0\to 3})} g(\eta_{0\to 2}) \cdot \ell_B (\eta_{2\to 3})\\
        &\overset{\mathclap{\strut(*)}}= \sum_{\substack{(\eta_{0\to 1}, \eta_{1\to 2}, \eta_{2\to 3})\in \textup{Decomp}^3_B(\eta_{0\to 3}):\\ \eta_{0\to 1} \textrm{ is algebraic, and }\eta_{1\to 2}\textrm{ is free}}}
        g(\eta_{0\to 2})\cdot  \ell_B (\eta_{2\to 3})\\
        &\overset{\mathclap{\strut(**)}}= \sum_{\substack{(\eta_{0\to 1}, \eta_{1\to 3})\in \textup{Decomp}_B(\eta_{0\to 3}):\\ \eta_{0\to 1} \textrm{ is algebraic}}} g(\eta_{0\to 1}) \sum_{\substack{(\eta_{1\to 2}, \eta_{2\to 3})\in \textup{Decomp}_B(\eta_{1\to 3}):\\ \eta_{1\to 2}\textrm{ is free} }} \ell_B (\eta_{2\to 3})\\
    \end{split}
\end{equation}
where in $(*)$ we condition on the algebraic-free decomposition (Theorem~\ref{thm_AFD}) of $\eta_{0\to 2}$ into its algebraic part $\eta_{0\to 1}$ and free part $\eta_{1\to 2}$, and in $(**)$ we use the free-invariance of $g$.
Using this, we prove the claims separately:
\begin{enumerate}
    \item By definition, if $\eta_{0\to 3}$ is algebraic then $\ell_{\textup{alg}}$ is the unique function satisfying
    \[ f(\eta_{0\to 3}) = \sum_{(\eta_{0\to 1}, \eta_{1\to 3})\in \textup{Decomp}_{\textup{alg}}(\eta_{0\to 3})} \ell_{\textup{alg}}(\eta_{1\to 3}) 
    = \sum_{\substack{(\eta_{0\to 1}, \eta_{1\to 3})\in \textup{Decomp}_B(\eta_{0\to 3}):\\ \eta_{0\to 1} \textrm{ is algebraic}}} \ell_{\textup{alg}}(\eta_{1\to 3}) \]
    where the last equality follows as $\eta_{1\to 3}$ must be algebraic too, by Observation~\ref{observe_alge_is_terminal_property}.
    Since the constant function $1$ is obviously free-invariant, we can substitute it for $g$ in equation~\eqref{eq_proof_of_lemma_mobius_alg} and get
    \[ f(\eta_{0\to 3}) = \sum_{\substack{(\eta_{0\to 1}, \eta_{1\to 3})\in \textup{Decomp}_B(\eta_{0\to 3}):\\ \eta_{0\to 1} \textrm{ is algebraic}}} \,\,\,\, \sum_{\substack{(\eta_{1\to 2}, \eta_{2\to 3})\in \textup{Decomp}_B(\eta_{1\to 3}):\\ \eta_{1\to 2}\textrm{ is free} }} \ell_B (\eta_{2\to 3})\]
    which implies
    \[ \ell_{\textup{alg}} (\eta_{1\to 3}) = \sum_{\substack{(\eta_{1\to 2}, \eta_{2\to 3})\in \textup{Decomp}_B(\eta_{1\to 3}):\\ \eta_{1\to2} \textrm{ is free}}} \ell_B (\eta_{2\to3}). \]
    \item Returning back to the more general free-invariant $g$, we now have
    \begin{equation*}
        \begin{split}
            (g\Bconv\mu^B\Bconv f)(\eta_{0\to 3})
            &= \sum_{\substack{(\eta_{0\to 1}, \eta_{1\to 3})\in \textup{Decomp}_B(\eta_{0\to 3}):\\ \eta_{0\to 1} \textrm{ is algebraic}}} g(\eta_{0\to 1}) \ell_{\textup{alg}}(\eta_{1\to 3})\\
            &= (g\convalg \ell_{\textup{alg}})(\eta_{0\to 3})\\
            &= (g \convalg \mu^{\textup{alg}} \convalg f)(\eta_{0\to 3}).
        \end{split}
    \end{equation*}
\end{enumerate}
\end{proof}

\begin{corollary}
(Algebraic Induction-Convolution Lemma)
\label{lemma_AlgICL}
Let $w\in F_r$ be a word, $G$ a compact group, and $\phi\in \textup{char}(G)$ a character. 
Assume that $\inner*{w}\le F_r$ is an algebraic extension.
Then
\[ \EX_{w\to F_r}[\textup{Ind}_n\phi] 
= \sum_{(\eta_1, \eta_2)\in \textup{Decomp}_{\textup{alg}}(w\to F_r)} \EX_{\eta_1}[\phi] \cdot L_{\eta_2}^{\textup{alg}}(\Sigma_n\acts[n]). \]
In the language of operators,
\begin{equation*}
    \EX\left[\textup{Ind}_{n}\phi\right] = \EX\left[\phi\right] \convalg \mu^{\textup{alg}} \convalg \EX\left[\Sigma_n\acts[n]\right]. 
\end{equation*}
\end{corollary}

\begin{proof}
    In Lemma~\ref{lemma_ICL} we proved 
    $ \EX[\textup{Ind}_n\phi] = \EX[\phi] \underset{B}{*} \mu^B \underset{B}{*} \EX[\textup{Ind}_n\textbf{1}] $,
    and in Lemma~\ref{lemma_mobius_alg_vs_B} (part 2) we proved that for every free-invariant $g\in \C^{\textup{sur}_B}$ and algebraic morphism $\eta$,
    $\left(g \Bconv \mu^B \Bconv f\right)(\eta) = \left(g \convalg \mu^{\textup{alg}} \convalg f\right)(\eta). $
    As $\EX\left[\phi\right]$ is free-invariant and $\eta_{w\to F_r}$ is algebraic, the desired result follows by substitution of $g = \EX\left[\phi\right], f = \EX[\textup{Ind}_n\textbf{1}]$ and $\eta = \eta_{w\to F_r}$.
\end{proof}

A special case of the algebraic induction-convolution lemma, where $G$ is cyclic and $\Sigma_n = S_n$, is proven in \cite[equation (3.7)]{magee2021surface}, where $L^{\textup{alg}}$ is called "contrib".

\section{Proof of Theorems \ref{thm_main_result}, \ref{thm_main_iterated_result}}
\label{section_proof_theorems}

\subsection{Proof of the Main Result}
\label{subsection_thm_main_result}

In this subsection we prove our main result, Theorem~\ref{thm_main_result}. It follows quite immediately from the induction-convolution lemma. 
Recall that for every compact group $G$ with an irreducible character $\phi\in \hat{G}$, we have the corresponding irreducible character of $G\wr S_n$,
\begin{equation*}
    \chi_{\phi, n} =
    \begin{cases}
    \textup{Ind}_n \phi & \textrm{if } \phi \neq \textbf{1}, \\
    \left(\textup{Ind}_n\textbf{1}\right) - \textbf{1} = \textup{std}\left(S_n\right) & \textrm{if } \phi = \textbf{1}.
    \end{cases}
\end{equation*}
This is an irreducible character of dimension $\phi(1)n - \bbone_{\phi=1}$.
Also recall that for a free word $w\in F_r$, we defined in Definition~\ref{def_witnesses_no_alg}
    \[\mathscr{C}_{\phi}(w) \defeq \sum_{H\in \textup{Crit}_{\phi}(w)} \EX_{w\to H}[\phi]. \]
In the following 2 propositions, we prove the main result (Theorem~\ref{thm_main_result}):

\begin{proposition}
\label{prop_rational_function}
(Rationality and Lower Bound)
Let $w\in F_r$ be a word, $G$ a compact group, and $\phi\in \hat{G}$. 
Then for $n\ge |w|$, $\EX_w[\textup{Ind}_n\phi]$ coincides with some rational function in $\Q_G(n)$. 
Moreover, if $\EX_w[\textup{Ind}_n\phi]\neq 0$ then 
\[\EX_w[\textup{Ind}_n\phi] = \Omega\left(n^{-|w|}\right). \]
\end{proposition}

\begin{proof}
Fix a basis $B\subseteq F_r$ such that $|w|_B = |w|$. 
By the induction-convolution lemma (Lemma~\ref{eq_ICL_withL}),
\[\EX_{w}[\textup{Ind}_n\phi] = \sum_{(\eta_1, \eta_2)\in \textup{Decomp}_B(w\to F_r)} \EX_{\eta_1}[\phi]\cdot L_{\eta_2}^B(n).\]
Now $\EX_{\eta_1}[\phi]\in \Q_G$, and by Proposition~\ref{prop_LB} if $n\ge |w|$ then $L_{\eta}^B(n)$ coincides with some rational function in $\Q(n)$. 

Now recall that for every rational function $P/Q$, where $P, Q$ are polynomials with $\deg(P)\le \deg(Q) = |w|$, the Laurent series $\frac{P(x)}{Q(x)} = \sum_{k=0}^{\infty} a_k x^{-k}$ satisfies a linear recurrence: if $Q(x) = \sum_{i=0}^{|w|} q_i x^i$ then for every $1\le k$, $\sum_{i=0}^{|w|} q_i a_{i + k} = 0$.
In particular, if $P$ is not the zero polynomial, then there exists $t\in \{1, \ldots, |w|\}$ such that $a_t\neq 0$, hence $\frac{P(x)}{Q(x)} = \Omega\left(x^{-t}\right) = \Omega\left(x^{-|w|}\right)$.
The second claim follows, since $\EX_{w}[\textup{Ind}_n\phi] = \frac{P(n)}{Q(n)}$ for some polynomials with degrees $\deg(P)\le \deg(Q) = |w|$.

\end{proof}

\begin{proposition}
(Approximation for Characters of Linear Dimension)\\
\label{prop_approx_lin_dim_char}
For every $w\in F_r$ and $\phi\in \hat{G}$,
\[ \EX_w[\chi_{\phi, n}] = \mathscr{C}_{\phi}(w)\cdot n^{1-\pi_{\phi}(w)} + O\left( n^{-\pi_{\phi}(w)} \right). \]
\end{proposition}

\begin{proof}
For the trivial character $\phi = \textbf{1}$ it follows from the main result of \cite{PP15}.
For every other $\phi$, the proof is very similar to the proof of \cite[Theorem 1.11]{magee2021surface}.
Consider the free-algebraic decomposition of $\eta_{w\to F_r}$: there is a subgroup $M$ satisfying $\inner*{w}\lealg M \leff F_r$.
Since $\EX\left[\textup{Ind}_{n}\phi\right]$ is free-invariant, $\EX_{w\to M}\left[\textup{Ind}_{n}\phi\right] = \EX_{w\to F_r}\left[\textup{Ind}_{n}\phi\right]$.
Thus we can assume without loss of generality that $\eta_{w\to F_r}$ is algebraic.
We apply the algebraic induction-convolution lemma (Lemma~\ref{lemma_AlgICL}) and get
\begin{equation*}
    \begin{split}
        \EX_w[\chi_{\phi, n}] = \EX_{w\to F_r}[\textup{Ind}_n\phi] 
        &= \sum_{(\eta_1, \eta_2)\in \textup{Decomp}_{\textup{alg}}(\eta_{w\to F_r})} \EX_{\eta_1}[\phi]\cdot L_{\eta_2}^{\textup{alg}}(n)\\
        &= \sum_{\substack{H:\quad w\in H\le F_r,\\ w\to H \textrm{ is algebraic}}} \EX_{w\to H}[\phi]\cdot L_{H\to F_r}^{\textup{alg}}(n).
    \end{split}
\end{equation*}
We care only about non-zero summands, that is, subgroups $H$ such that $\EX_{w\to H}[\phi]\neq 0$.
The subgroup $H = \inner*{w}$ is absent among these subgroups: $\EX_{w\to w}[\phi] = \inner*{\phi, 1} = 0$, since $\phi$ is irreducible and non trivial.
Thus we sum only over subgroups $H$ that are proper algebraic extensions of $w$ that also satisfy $\EX_{w\to H}[\phi]\neq 0$: in other words, elements of $\textup{Wit}_{\phi}(w)$. 
Among witnesses subgroups, we are interested only in those with minimal rank, as by Corollary~\ref{corollary_L_alg_approx},
$ L_{H\to F_r}^{\textup{alg}}(n) = n^{1-\textup{rk}(H)} \cdot \left(1 + O\left(n^{-1}\right)\right), $
so we may restrict the sum further to $H\in \textup{Crit}_{\phi}(w)$ only, in which every subgroup has Euler characteristic $1 - \pi_{\phi}(w)$, and get
\begin{equation*}
    \begin{split}
        \EX_w[\chi_{\phi, n}] 
        &= \sum_{H\in \textup{Crit}_{\phi}(w)} \EX_{w\to H}[\phi]\cdot n^{1 - \pi_{\phi}(w)} \cdot \left(1 + O\left(n^{-1}\right)\right)\\
        (\textup{Definition}~\ref{def_witnesses_no_alg})\quad\quad\quad\quad &= \mathscr{C}_{\phi}(w)\cdot n^{1-\pi_{\phi}(w)} + O\left( n^{-\pi_{\phi}(w)} \right)
    \end{split}
\end{equation*}
which implies the desired result for $\phi \neq \textbf{1}$. 
\end{proof}

We do not know what are the possible values of $\pi_{\phi}$, except for linear characters $\phi$:
\begin{remark}
\label{remark_composition_with_balanced_morphisms}
Let $m\in \{2, 3, \ldots\}\sqcup \{\infty\}$, let $\phi\colon C_m\to \Ss^1$ be the standard embedding, and let $w\in H\le J\le F_r$ be a chain of $B$-surjective extensions (for some basis $B\subseteq F_r$). 
Recall $K_m$ from Definition~\ref{def_K_m(H)}. 
It is clear that $K_m(H)\le K_m(J)$, and together with \cite[Lemma 3.2]{magee2021surface} and the example in equation~\eqref{equation_expectation_of_cyclic_embedding}, we get that if $\EX_{\eta_{w\to H}}[\phi]\neq 0$, then also $\EX_{\eta_{w\to J}}[\phi]\neq 0$. 
As a consequence, for every $w\in F_r$, $\pi_{C_m}(w)\in \{0, 1, \ldots, r\} \sqcup \{\infty\}$.
Indeed, every $H\in \textup{Crit}_{C_m}(w)$ has a $B$-surjective map to $F_r$, so $\EX_{\eta_{w\to F_r}}[\phi]\neq 0$ whenever $\textup{Crit}_{C_m}(w)\neq \emptyset$.

In contrast, for non-linear characters $\phi$, it is possible that $\EX_{\eta_{w\to H}}[\phi]\neq 0$ but also $\EX_{\eta_{w\to J}}[\phi] = 0$:
Let $n\ge 3$, let $\phi = \#\textup{fix}-1 \in \widehat{S_n}$ be the standard character of $S_n$, and denote $\ell\defeq \textup{lcm}(1, \ldots, n)$.
We define a word $w\in F_2 = \textup{Free}(\{x, y\})$ by $w = x^{-3}\left(xy^{\ell}\right)^2$, and a subgroup $H\defeq \inner*{x, y^{\ell}} \le F_2$ containing $w$.
Then $\EX_{w\to H}[\phi] = \frac{1}{n-1}$, whereas $\EX_{w\to F_2}[\phi] = 0$. 
Indeed, by \cite{frobenius1896gruppencharaktere}, 
\[\EX_{w\to H}[\phi] 
= \frac{1}{\dim\phi} \EX_{\sigma\sim U(S_n)} [\phi(\sigma^{-3})] \EX_{\tau\sim U(S_n)}[ \phi (\tau^2)] 
= \frac{1}{n-1} (2-1)\cdot (2-1) = \frac{1}{n-1},\]
and since $\ell$ is the exponent of $S_n$, $\EX_{w\to F_2}[\phi] = \EX_{\sigma\sim U(S_n)} [\phi(\sigma^{-3+2})] = 1-1 = 0$.


\end{remark}

\subsection{Iterated Wreath Products}
\label{subsection_iterated_wreath}

Recall Definition~\ref{def_iterated_W_Ind}. Our goal now is to prove Theorem~\ref{thm_main_iterated_result}.
We have an iterated version of the induction-convolution lemma, suited for iterated wreath products.
Recall $\textup{Decomp}_{\textup{alg}}^m$ from Definition~\ref{def_decomp_alg_B}.

\begin{lemma}
\label{lemma_algebraic_iterated_induction_convolution}
(Iterated Algebraic Induction-Convolution Lemma)
Let $G$ be a compact group with a character $\phi\in \textup{char}(G)$. Then we have the following equality of operators:
\begin{equation*}
    \EX\left[\textup{Ind}_{n_1, \ldots, n_m}\phi\right] = \EX\left[\phi\right] \convalg \mu^{\textup{alg}} \convalg \EX\left[S_{n_1}\acts [n_1]\right] \convalg \ldots \convalg \mu^{\textup{alg}} \convalg \EX\left[S_{n_m}\acts [n_m]\right]. 
\end{equation*}
In plain words, let $w\in F_r$ be a word such that $F_r$ is an algebraic extension of $\inner*{w}$. Then
\begin{equation*}
    \EX_{w}\left[\textup{Ind}_{n_1, \ldots, n_m}\phi\right] = \sum_{\substack{(\eta_0, \ldots, \eta_m)\in \\\textup{Decomp}_{\textup{alg}}^{m + 1} \left(w\to F_r\right)}} \EX_{\eta_0}\left[\phi\right] \cdot L_{\eta_1}^{\textup{alg}}\left(n_1\right)\cdots L_{\eta_m}^{\textup{alg}}\left(n_m\right). 
\end{equation*}
\end{lemma}

\begin{proof}
It follows by iterating $m$ times the algebraic induction-convolution lemma (\ref{lemma_AlgICL}).
\end{proof}

We state again Theorem~\ref{thm_main_iterated_result} for convenience: let $G$ be a compact group, let $\phi\in \hat{G}-\{\textbf{1}\}$, let $w\in F_r$ and fix $m\in \N$. 
For every sequence $n_1, \ldots, n_m\in \N$ satisfying $n_i \geq |w|$ for $i=1,\ldots, m$,
denote $\chi\defeq \textup{Ind}_{n_1, \ldots, n_m}\phi \in \widehat{\mathcal{W}_{n_1, \ldots, n_m}(G)}$.
Then $\EX_w\left[\chi\right]$ coincides with some rational function in $\Q_G(n_1, \ldots, n_m)$, and 
\[ \EX_w\left[\chi\right] = O\left(\dim(\chi)^{1-\pi(w)}\right) \]
where the implied constant depends on $w, m$ and $G$.
Moreover, if $\phi$ is a linear character, then 
\[ \EX_w\left[\chi\right] = (n_1\cdots n_m)^{1-\pi_{\phi}(w)}\cdot \left(C + O\left(\sum_{i=1}^m n_i^{-1}\right)\right), \]
for some positive integer $C\le |\textup{Crit}_{\phi}(w)|^m$ which depends on $w, m$.


\begin{proof}
[Proof of Theorem~\ref{thm_main_iterated_result}]
The proof is very similar to the proof of the main result.
By Lemma~\ref{lemma_algebraic_iterated_induction_convolution}, 
\begin{equation*}
    \EX_{w}\left[\textup{Ind}_{n_1, \ldots, n_m}\phi\right] = \sum_{\substack{(\eta_0, \ldots, \eta_m)\in \\\textup{Decomp}_{\textup{alg}}^{m + 1} \left(w\to F_r\right)}} \EX_{\eta_0}\left[\phi\right] \cdot L_{\eta_1}^{\textup{alg}}\left(n_1\right)\cdots L_{\eta_m}^{\textup{alg}}\left(n_m\right),
\end{equation*}
which gives rationality as each $L_{\eta_i}^{\textup{alg}}\left(n_{i}\right)$ coincides with some rational function.

Among the chains of morphisms in $\textup{Decomp}_{\textup{alg}}^{m + 1} \left(w\to F_r\right)$, we are interested only in chains with $\EX_{\eta_0}\left[\phi\right]\neq 0$, so $\textup{Im}(\eta_0)\in \textup{Wit}_{\phi}(w)$ (as $w$ itself has no contribution, since $\inner*{\phi, \textbf{1}} = 0$). 
If $\phi$ is linear, then all the subgroups $\textup{Im}(\eta_i)\in \textup{Wit}_{\phi}(w)$ are $\phi$-witnesses (by Remark~\ref{remark_composition_with_balanced_morphisms}), and their Euler characteristics are at most $1-\pi_{\phi}(w)$; 
otherwise, even if $\phi$ is not linear, we still know that $\textup{Im}(\eta_i)$ are non-trivial algebraic extensions of $w$, so their Euler characteristics are at most $1-\pi(w)$.
In the first case, the desired constant $C$ is the number of chains in $\textup{Decomp}_{\textup{alg}}^{m + 1} \left(w\to F_r\right)$ that pass only through $\phi$-critical subgroups (since in this case $\EX_{\eta_0}\left[\phi\right] = 1$).
Since $L_{H\to J}^{\textup{alg}}\left(n_{i}\right) = n^{1-\textup{rk}(H)}\cdot \left(1 + O\left(n_i^{-1}\right)\right)$, the result follows.

\end{proof}

\subsubsection{Spherically Symmetric Trees}
\label{subsection_spherical_tree}

The procedure of iterating $\textup{Ind}_{n_i}$ yields characters that correspond to actions on trees, in the following sense.

\begin{definition}
(The Spherically Symmetric Tree $\mathcal{T}_{n_1, \ldots, n_m}$) Let $n_1, \ldots, n_m\in \N$. The \textbf{spherically symmetric tree} $\mathcal{T}_{n_1, \ldots, n_m}$ is a rooted tree with $m + 1$ layers, indexed by $0, \ldots, m$.
The $0$-layer contains a single vertex - the root. 
For every $1\le i\le m$, every vertex in the $(i-1)$-layer has $n_{m-i}$ children in the $i$-layer, numbered from $1$ to $n_{m-i}$.
Thus the $i$-layer can be identified with $[n_m]\times [n_{m-1}]\times \ldots \times [n_{m-i+1}]$. 
In particular there are $n_1\cdots n_m$ leaves (in the $m$-layer).
\end{definition}

\begin{example}
The tree $\mathcal{T}_{3, 2}$:
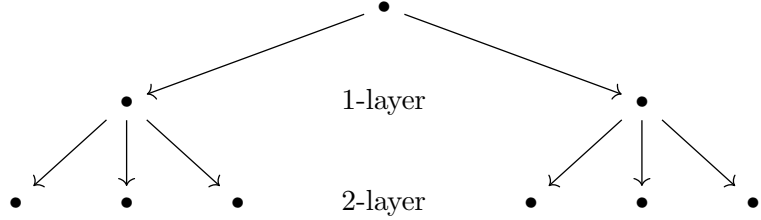
\[\begin{tikzcd}
	&&& \bullet \\
	& \bullet && {1\textup{-layer}} && \bullet \\
	\bullet & \bullet & \bullet & {2\textup{-layer}} & \bullet & \bullet & \bullet
	\arrow[from=1-4, to=2-2]
	\arrow[from=1-4, to=2-6]
	\arrow[from=2-2, to=3-1]
	\arrow[from=2-2, to=3-2]
	\arrow[from=2-2, to=3-3]
	\arrow[from=2-6, to=3-5]
	\arrow[from=2-6, to=3-6]
	\arrow[from=2-6, to=3-7]
\end{tikzcd}\]
\end{example}

\begin{observation}
The group of automorphisms of $\mathcal{T}_{n_1, \ldots, n_m}$ is the iterated wreath product $\mathcal{W}_{n_1, \ldots, n_m}\defeq \mathcal{W}_{n_1, \ldots, n_m}(\{1\}) = \left(\ldots\left(S_{n_1}\wr S_{n_2}\right)\wr S_{n_3}\ldots \right)\wr S_{n_m}$. 
\end{observation}

Note that $\left|\mathcal{W}_{n_1, \ldots, n_m}\right| = (n_1!)^{n_2\cdots n_m}\cdot (n_2!)^{n_3\cdots n_m}\cdots n_m!.$
An element of $\mathcal{W}_{n_1, \ldots, n_m}$ can be described as a function $\vec{f}\colon V(\mathcal{T})-\{\textup{leaves}\}\to \bigsqcup_{i=1}^{m} S_{n_i}$, sending the vertices of the $i$-layer to $S_{n_{m-i}}$, for $0\le i \le m-1$.
Explicitly, 
$\vec{f} = (f_1, \ldots, f_{m})$, where for every $1\le i \le m$, $f_i\colon [n_m]\times [n_{m-1}]\times \ldots \times [n_{i+1}] \to S_{n_i}$, which can be visualized as putting an $S_{n_i}$-element on every vertex in the $(m-i)$-layer of $\mathcal{T}_{n_1, \ldots, n_m}$. In particular, $f_{m}\in S_{n_{m}}$ as there is only one vertex in the $0$-layer.
In the example above with $m = 2, n_1 = 3, n_2 = 2$, a tree automorphism is $(f_1, f_2)$ with $f_1\colon [2]\to S_3$ and $f_2\in S_2$.

What is the decomposition of the permutation representation $\mathcal{W}_{n_1, \ldots, n_m}\acts \mathcal{T}_{n_1, \ldots, n_m}$?
The orbits are the layers, so we may restrict attention to the action on the leaves (the last layer), or equivalently $\mathcal{W}_{n_1, \ldots, n_m}\acts [n_1]\times \ldots \times [n_m]$. 

\begin{proposition}
The character of the permutation representation $\mathcal{W}_{n_1, \ldots, n_m}\acts [n_1]\times \ldots \times [n_m]$ is $\textup{Ind}_{n_1, \ldots, n_m}\textbf{1}$, where $\textbf{1}$ is the trivial character of the trivial group. 
\end{proposition}

\begin{proof}
By induction on $m$. The case $m=0$ is trivial.
The number of fixed points of an element $\vec{f}\in \mathcal{W}_{n_1, \ldots, n_m}$ is $\sum_{i\in [n_m]: f_m.i = i} \#\textup{fix}\left(\vec{f}\restriction_{i}\right)$, where $\#\textup{fix}\left(\vec{f}\restriction_{i}\right)$ is the number of fixed points of the restriction of $\vec{f}$ to the sub-tree of the $i^{th}$ child of the root. 
The restricted function $\vec{f}\restriction_{i}$ can be naturally considered as an element of $\mathcal{W}_{n_1, \ldots, n_{m-1}}(\{1\})$ acting on this $i^{th}$ sub-tree, so by the induction hypothesis, $\#\textup{fix}\left(\vec{f}\restriction_{i}\right) = \textup{Ind}_{n_1, \ldots, n_{m-1}}\textbf{1}(\vec{f}\restriction_{i})$, and by definition
$\sum_{i\in [n_m]: f_m.i = i} \#\textup{fix}\left(\vec{f}\restriction_{i}\right) = \textup{Ind}_{n_m} \#\textup{fix}\left(\vec{f}\right)$, concluding the proof.
\end{proof}

\begin{observation}
\label{observe_decompose_Ind1}
The decomposition of $\textup{Ind}_{n_1, \ldots, n_m}\textbf{1}$ to irreducible characters is
\[\textup{Ind}_{n_1, \ldots, n_m}\textbf{1} = \textbf{1} + \textup{std}_{n_m} + \textup{Ind}_{n_m}\textup{std}_{n_{m-1}} + \ldots + \textup{Ind}_{n_2, \ldots, n_m}\textup{std}_{n_1},\]
where $\textup{std}_{n_i}\in \widehat{S_{n_i}} = \#\textup{fix} - 1$ is the standard $(n_i-1)$-dimensional irreducible character, and for every $1\le i\le m$ the character $\textup{Ind}_{n_{i + 1}, \ldots, n_m}\textup{std}_{n_i}$ is inflated to $\mathcal{W}_{n_1, \ldots, n_m}$ by the epimorphism $\mathcal{W}_{n_1, \ldots, n_m} \twoheadrightarrow \mathcal{W}_{n_{i+1}, \ldots, n_m}$.
Note that the dimensions fit:
\[n_1\cdots n_m = 1 + (n_m - 1) + n_m(n_{m-1} - 1) + \ldots + n_m\cdots n_2(n_1 - 1).\]
\end{observation}

\begin{proof}
By Proposition~\ref{prop_Ind_phi_is_irred}, the operator \textup{Ind} gives the following branching diagram of irreducible characters:

\[\begin{tikzcd}
	{\widehat{\mathcal{W}_{n_1, \ldots, n_m}}:} & {\textbf{1}} & {\textup{std}_{n_m}} & {\textup{Ind}_{n_3, \ldots, n_m}\textup{std}_{n_2}} & {\textup{Ind}_{n_2, \ldots, n_m}\textup{std}_{n_1}} \\
	{\widehat{S_{n_1}\wr S_{n_2}}:} & {\textbf{1}} & {\textup{std}_{n_2}} & {\textup{Ind}_{n_2}\textup{std}_{n_1}} \\
	{\widehat{S_{n_1}}:} & {\textbf{1}} & {\textup{std}_{n_1}} \\
	{\widehat{\{1\}}:} & {\textbf{1}}
	\arrow[from=4-2, to=3-2]
	\arrow[from=4-2, to=3-3]
	\arrow[from=3-2, to=2-2]
	\arrow[from=3-2, to=2-3]
	\arrow["\vdots"{description}, draw=none, from=2-2, to=1-2]
	\arrow[from=3-3, to=2-4]
	\arrow["\cdots"{marking}, draw=none, from=2-2, to=1-3]
	\arrow["\cdots"{marking}, draw=none, from=2-4, to=1-5]
	\arrow["\cdots"{marking}, draw=none, from=2-3, to=1-4]
	\arrow["\cdots"{marking}, draw=none, from=1-3, to=1-4]
	\arrow["\vdots"{description}, draw=none, from=2-1, to=1-1]
	\arrow["{\textup{Ind}_{n_1}}", from=4-1, to=3-1]
	\arrow["{\textup{Ind}_{n_2}}", from=3-1, to=2-1]
\end{tikzcd}\]
\end{proof}

Now we can prove Corollary~\ref{corollary_intro_spherical_tree}: for every $w\in F_r$,
\[ \left| \EX_w\left[\textup{Ind}_{n_1, \ldots, n_m}\textbf{1}\right] - \EX_w\left[\#\{\textup{fix}\left(S_{n_m}\acts [n_m]\right)\}\right] \right| = O\left((n_{m-1} n_m)^{1-\pi(w)}\right)_{n_m\to\infty}. \]

\begin{proof}
[Proof of Corollary~\ref{corollary_intro_spherical_tree}]
Compare the decomposition of $\textup{Ind}_{n_1, \ldots, n_m}\textbf{1}$ into irreducible character to the decomposition of the fixed points counting character of $S_{n_m}$: the difference consist only of characters with dimension $\ge n_m \cdot n_{m-1}$ and small $w$-measure.
\begin{equation*}
    \begin{split}
        \EX_w\left[\textup{Ind}_{n_1, \ldots, n_m}\textbf{1} - S_{n_m}\acts [n_m]\right] &= \\
        (\textup{Observation}~\ref{observe_decompose_Ind1})\quad\quad &= \sum_{i=1}^{m-1} \EX_w\left[ \textup{Ind}_{n_{i + 1}, \ldots, n_m}\textup{std}_{n_i} \right] \\
        (\textup{Theorem}~\ref{thm_main_iterated_result}) \quad\quad&= O\left((n_{m-1} n_m)^{1-\pi(w)}\right).
    \end{split}
\end{equation*}
\end{proof}

\section{Computation of Witnesses}
\label{section_witnesses}

In this chapter we elaborate more about witnesses, and prove Theorem~\ref{thm_pi_phi_inequalities}. 

\begin{proposition}
\label{prop_pi_phi_of_independent_words}
    Let $G$ be a compact group and let $\phi\in \hat{G}-\{\textbf{1}\}$.
    Let $w_1, w_2\in F_r$ be "disjoint" words: that is, there is a free decomposition $F_r = J_1 * J_2$ such that $w_i$ can be conjugated into $J_i$, for $i=1, 2$.
    Then
        \begin{align*}
        \pi_{\phi}(w_1w_2) & = \pi_{\phi}(w_1) + \pi_{\phi}(w_2), \\
        \mathscr{C}_{\phi}(w_1 w_2) & = \frac{\mathscr{C}_{\phi}(w_1) \mathscr{C}_{\phi}(w_2)}{\dim\phi}, \\
        \textup{Crit}_{\phi}(w_1 w_2) & = \{H_1 * H_2: \,\,H_i\in \textup{Crit}_{\phi}(w_i)\}. 
    \end{align*}

\end{proposition}

\noindent The case $\phi = \textbf{1}$ is stated already in \cite[Lemma 6.8]{Puder_2014}, claiming both
\[ \pi(w_1 w_2) = \pi(w_1) + \pi(w_2)\quad \textup{ and } \quad \textup{Crit}(w_1 w_2) = \{H_1 * H_2: H_i\in \textup{Crit}(w_i)\}. \]
Neither the proof nor its main tool, which is the following proposition, has been published.
\begin{proposition}
\label{prop_alg_disjoint_words}
    (D. Puder, personal communication, March 25, 2023)
    Let $w_1, w_2$ be disjoint words.
    Denote by $\textup{Alg}(w)$ the set of all algebraic extensions of $w$.
    Then the free product map $*\colon \textup{Alg}(w_1)\times \textup{Alg}(w_2)\to \textup{Alg}(w_1 w_2)-\{\inner*{w_1 w_2}\}$ is bijective.
\end{proposition}

\begin{proof} [Proof of Proposition~\ref{prop_pi_phi_of_independent_words}]
By \cite{frobenius1896gruppencharaktere} (see also \cite{magee2021surface}), if $H\le F_r$ contains both $w_1, w_2$ then
\begin{equation}
\label{eq_frobenius_indep_words_formula}
   \EX_{w_1w_2\to H}[\phi] = \frac{\EX_{w_1\to H}[\phi]\EX_{w_2\to H}[\phi]}{\dim\phi}. 
\end{equation}
As a consequence, $\EX_{w_1w_2\to H}[\phi]\neq 0$ if and only if $\EX_{w_i\to H}[\phi] \neq 0$ for both $i = 1, 2$.
Moreover, if $H = H_1 * H_2$ for $H_i$ containing $w_i$, then $\EX_{w_i\to H}[\phi] = \EX_{w_i\to H_i}[\phi]$.
This lets us restrict the free-product map from Proposition~\ref{prop_alg_disjoint_words} to $\phi$-witnesses only:
\[ * \colon \textup{Wit}_{\phi}(w_1)\times \textup{Wit}_{\phi}(w_2)\to \textup{Wit}_{\phi}(w_1w_2). \]
This restricted map is still bijective, as $\textup{Wit}_{\phi}(w_1w_2)$ does not contain $\inner*{w_1 w_2}$.
Since $\textup{rk}(H_1 * H_2) = \textup{rk}(H_1) + \textup{rk}(H_2)$, the subgroups of minimal rank in $\textup{Wit}_{\phi}(w_1w_2)$ are precisely the free products of subgroups of minimal rank, so we can restrict the map to critical subgroups and get the bijective map 
\begin{equation}
\label{eq_crit_disjoint_words}
    * \colon \textup{Crit}_{\phi}(w_1)\times \textup{Crit}_{\phi}(w_2)\to \textup{Crit}_{\phi}(w_1w_2).
\end{equation}
In particular we get $\pi_{\phi}(w_1w_2) = \pi_{\phi}(w_1) + \pi_{\phi}(w_2)$.
It is left to compute $\mathscr{C}_{\phi}(w_1 w_2)$:
\begin{equation*}
    \begin{split}
        \mathscr{C}_{\phi}(w_1 w_2)
        &= \sum_{H\in \textup{Crit}_{\phi}(w_1 w_2)} \EX_{w_1 w_2\to H}[\phi] \\
        (\textup{equation}~\eqref{eq_crit_disjoint_words}) &= \sum_{H_1\in \textup{Crit}_{\phi}(w_1)} \sum_{H_2\in \textup{Crit}_{\phi}(w_2)} \EX_{w_1 w_2\to H_1 * H_2}[\phi] \\
        (\textup{equation}~\eqref{eq_frobenius_indep_words_formula}) &= \sum_{H_1\in \textup{Crit}_{\phi}(w_1)} \sum_{H_2\in \textup{Crit}_{\phi}(w_2)} \frac{\EX_{w_1\to H_1 * H_2}[\phi]\EX_{w_2\to H_1 * H_2}[\phi]}{\dim\phi} \\
        (\textup{Definition}~\ref{def_free-inv_func}) &= \frac{1}{\dim\phi} \left(\sum_{H_1\in \textup{Crit}_{\phi}(w_1)} \EX_{w_1\to H_1}[\phi]\right)\left(\sum_{H_2\in \textup{Crit}_{\phi}(w_2)} \EX_{w_2\to H_2}[\phi]\right) \\
        &= \frac{\mathscr{C}_{\phi}^{\pi}(w_1)\mathscr{C}_{\phi}^{\pi}(w_2)}{\dim\phi} . \\
    \end{split}
\end{equation*}

\end{proof}

It is natural to ask: when is the inequality $\pi(w)\le \pi_{\phi}(w)$ tight? 
Apparently, for every $w$ there are infinitely many characters $\phi$ with $\pi(w) = \pi_{\phi}(w)$; specifically, standard characters $\textup{std}\left(S_n\right)$ of symmetric groups:

\begin{proposition}
For every $w\in F_r$, the set 
$\{n\in \N: \pi(w) = \pi_{\textup{std}\left(S_n\right)}(w) \} $
is co-finite. 
The minimum of this set is at most $|w|$, and in fact at most
$ \max\{\textup{len}_H(w)\mid H\in \textup{Crit}(w)\} $
where $\textup{len}_H(w)$ is the minimal length of $w$ with respect to bases of $H$.
(Here $\textup{std}\in \hat{S}_n$ is the standard representation $\#\left\{\textrm{fixed points}\right\} - 1$).
\end{proposition}

\begin{proof}
If $w\in F_r$ is primitive then we always have equality $\pi(w) = \pi_{\textup{std}\left(S_n\right)}(w) = \infty$, so we may assume that there exists a critical group $w\in H$ of rank $\pi(w)$ in which $w$ is not primitive.
By \cite{PP15}, the property "$w\in H$ is not primitive" is equivalent to both
\begin{itemize}
    \item $\{n\in \N: \EX_{w\to H}[\textup{std}\left(S_n\right)] \neq 0\}$ is co-finite (\cite[Theorem (1.4’)]{PP15}),
    \item $\min\{n\in \N: \EX_{w\to H}[\textup{std}\left(S_n\right)] \neq 0\} \le \textup{len}_H(w)$ (\cite[Proposition 1.6]{PP15}).
\end{itemize}
Now let $n\in \N$. If $\EX_{w\to H}[\textup{std}\left(S_n\right)] \neq 0$, then $\pi(w) \le \pi_{\textup{std}\left(S_n\right)}(w)\le \textup{rk}(H) = \pi(w). $

To complete the proof, $\textup{len}_H(w)\le \textup{len}(w)$ for every $w$-critical group $H$ since the morphism $w\to H$ is $B$-surjective so every letter of $w$ can increase its $H$-length by at most $1$.

\end{proof}

\begin{corollary}
\label{corollary_pi_limit_std}
For every $w\in F_r$,
\[\pi(w) = \lim_{n\to \infty} \pi_{\textup{std}\left(S_n\right)}(w) = \min_{n\in \N} \pi_{\textup{std}\left(S_n\right)}(w). \]
\end{corollary}
\noindent (Note that we do not claim that there is $n\in \N$ such that $\pi = \pi_{\textup{std}\left(S_n\right)}$).

\noindent Next, we show a bound on $\pi_{\phi}$ for characters of nilpotent groups, using a result of \cite{COCKE2019440}:

\begin{proposition}
Let $G$ be a finite nilpotent group, with exponent $\exp(G)$.
Then for every word $w\in F_r$:
\begin{equation*}
    \bigcup_{\phi\in \hat{G}-\{\textbf{1}\}} \textup{Wit}_{\phi}(w) = \bigcup_{p\divides exp(G) \textrm{ prime}} \textup{Wit}_{C_p}(w).
\end{equation*}
In particular, by taking minimal ranks,
\[ \min\{\pi_{\phi}(w) : \phi\in \hat{G}-\{\textbf{1}\}\} = \min\{\pi_{C_p}(w): p \textup{ is a prime divisor of }\exp(G) \}. \]
(One can replace $\exp(G)$ by $|G|$).
\end{proposition}

\begin{proof}
For every subgroup $H\le F_r$ define a function $\gcd_H \colon H\to \Z$ by composing the quotient homomorphism $H\twoheadrightarrow H/[H, H] \cong \Z^{\textup{rk}(H)}$ with the function $\textup{gcd}\colon \Z^{\textup{rk}(H)}\to \Z$ that computes the greatest common divisor of ${\textup{rk}(H)}$ integers.
Also denote $\gcd(w) = \gcd_{F_r}(w)$.
We cite a special case of \cite[Lemma 7]{COCKE2019440}:
let $G$ be a finite nilpotent group, and let $w\in H$ be a word in a f.g.\ free group. 
The following are equivalent:
\begin{enumerate}
    \item The word measure $\mu_{w\to H}$ induced by $w$ on $G$ is uniform.
    \item $\gcd(\gcd_H(w), \exp(G)) = 1$.
\end{enumerate}
To this end, fix $w$ and consider only proper algebraic extensions $H$ of $\inner*{w}$.
Note that the first property is equivalent to 
$\forall \phi\in \hat{G}-\{\textbf{1}\}:\,\, \EX_{w\to H}[\phi] = 0$, 
i.e.\ $H\not\in \textup{Wit}_{\phi}(w)$, and the last property is equivalent to $\forall p \divides \exp(G)$ prime, $\gcd_H(w)\not \equiv _p 0$, i.e.\ $w\not \in K_p(H)$.
So if the primes that divide $\exp(G)$ are $p_1, \ldots, p_k$, we have
\begin{equation*}
    \begin{split}
        \bigcup_{\phi\in \hat{G}-\{\textbf{1}\}} \textup{Wit}_{\phi}(w) 
        &= \{H \textup{ proper algebraic extension of } w \textup{ such that } \exists i: w\in K_{p_i}(H)\}\\
        &= \{H \textup{ proper algebraic extension of } w \textup{ such that } \exists i: H\in \textup{Wit}_{C_{p_i}}(w)\}.
    \end{split}
\end{equation*}
This finishes the proof. \\
\end{proof}

To deduce a simple corollary, let $G$ be a finite $p$-group (for some prime $p$). It is of course nilpotent. 
It should not surprise us that there is a character $\phi\in \hat{G}$ with $\pi_{\phi} = \pi_{C_p}$; indeed, $\Z/p\Z$ is a quotient of $G$ (by induction on the group size). 
But the previous theorem gives us the following non-obvious inequality:
\begin{corollary}
\label{corollary_pi_phi_of_nilpotent}
For every $\phi\in \hat{G}-\{\textbf{1}\}$,
$\pi_{\phi}\ge \pi_{C_p}. $
\end{corollary}

\appendix

\section{General Wreath Products}
\label{section_general_wreath}

In this appendix we demonstrate \textbf{asymptotically oligomorphic} sequences of group actions $\left(\Sigma_n\acts X_n\right)_{n=1}^{\infty}$, and show that they may replace the sequence $S_n\acts [n]$ in Theorem~\ref{thm_main_result} at the cost of weakening the bound on $\EX_w\left[ \textup{Ind}_{X_n}\phi \right]$.
Given a group action $\Sigma\acts X$ we denote its set of orbits by $X/\Sigma$.

\subsection{Asymptotically Oligomorphic Sequences}

 There is a known concept of \textbf{oligomorphic} group actions: a group action $\Sigma\acts X$ is called oligomorphic if the diagonal action $\Sigma\acts X^t$ has only finitely many orbits for every $t\in \N$. See \cite{cameron_1990} for a survey.
We define a very similar notion:

\begin{definition}
\label{def_asymp_oligom}
A sequence of finite group actions $\{\Sigma_n\acts X_n\}_{n=1}^{\infty}$ is called \textbf{asymptotically oligomorphic} if for every $t\in \N$, the sequence of diagonal actions $\Sigma_n\acts X_n^t$ has $O(1)$ orbits:
\[\sup_{n} \left| X_n^t/\Sigma_n \right| < \infty. \]
\end{definition}

\begin{example}
\label{example_oligom}
Some examples and non-examples for asymptotically oligomorphic actions are:
\begin{enumerate}
    \item For every $k$, the action of the symmetric group $S_n$ on subsets of size $k$ of $\{1, \ldots, n\}$ is asymptotically oligomorphic.
    Indeed, by \cite[Claim 4.7]{reiter19}, the number of orbits of $S_n\acts \binom{[n]}{k}^t $ is at most
    $\left(k + 1\right)^{2^{t} - 1}$.
    
    \item For every finite field $\F_q$ (where $q$ is some prime power indicating the field size), the sequence of actions $\left(\textup{GL}_n(\F_q)\acts \F_q\right)_{n=1}^{\infty}$ is asymptotically oligomorphic. In fact, for every $m\in \N$, the sequence of actions $\left(\textup{GL}_n\left(\F_q\right)\acts M_{n\times m}(\F_q)\right)_{n=1}^{\infty}$ by left multiplication (where $M_{n\times m}(\F_q)$ is the set of all $n\times m$ matrices) is also asymptotically oligomorphic.
    
    Indeed, for every $t\in \N$, consider the action $\textup{GL}_n\left(\F_q\right)\acts \left(M_{n\times m}(\F_q)\right)^t$: it is isomorphic to the action $\textup{GL}_n(\F_q)\acts M_{n\times mt}(\F_q)$, and if $n\ge mt$, then $B, C \in M_{n\times mt}$ are in the same orbit if and only if they have the same right kernel $\ker(B) = \ker(C)$, so the number of orbits is the number of subspaces of $\F_q^{mt}$: in particular, if $n\ge mt$ then the number of orbits does not depend on $n$. 
    (This number of subspaces is also known as $\sum_{k=0}^{mt} \genfrac{[}{]}{0pt}{}{mt}{k}_q$).
    \item The sequence $\left(\textup{GL}_n(\F_q)\acts \F_q^n\right)_{q\to \infty}$ is \textbf{not} asymptotically oligomorphic.
    Indeed, assume $n < t$ 
    and let $e_1, \ldots, e_n$ be a basis of $\F_q^n$. Then the function $(\F_q^n)^{t-n} \hookrightarrow (\F_q^n)^t, (v_{n+1}, \ldots, v_t) \mapsto (e_1, \ldots, e_n, v_{n+1}, \ldots, v_t)$ sends each vector into a different $\textup{GL}_n(\F_q)$-orbit, so there are at least $q^{n(t-n)}$ orbits. 
\end{enumerate}
\end{example}

Asymptotically oligomorphic sequences naturally arise as exhaustion of oligomorphic groups by finite subgroups:

\begin{observation}
Let $\Sigma_{\infty}\acts X_{\infty}$ be an oligomorphic group action, and let $\{\Sigma_n\}_{n\in\N}$ be an ascending sequence of finite subgroups (that is, $\Sigma_n\le \Sigma_{n+1}$) that exhaust $\Sigma_{\infty}$: $\bigcup_{n=1}^{\infty} \Sigma_n = \Sigma_{\infty}$. 
Then there is an exhaustion of $X_{\infty}$ by finite sets $\{X_n\}_{n=1}^{\infty}$ such that $\Sigma_n$ preserves $X_n$ and the sequence $\left(\Sigma_n\acts X_n\right)_{n=1}^{\infty}$ is asymptotically oligomorphic.

Indeed, let $\{x_i\}_{i=1}^{|X_{\infty} / \Sigma_{\infty}|}$ be a set of representatives of the orbits $X_{\infty} / \Sigma_{\infty}$, and define $X_n\defeq \bigcup_{i=1}^{|X_{\infty} / \Sigma_{\infty}|} \Sigma_n(x_i)$. Now for every $t\in \N$, there is a finite set of representatives $\mathscr{R}_t$ of $X_{\infty}^t / \Sigma_{\infty}$, and for every large enough $n$, $\mathscr{R}_t\subseteq X_n^t$ and $|X_n^t / \Sigma_n| = |\mathscr{R}_t|$ so $\sup_n |X_n^t / \Sigma_n| < \infty$.
\end{observation}

\begin{remark}
For every ascending sequence $\{X_n\}_{n=1}^{\infty}$ exhausting $X_{\infty}$ (that is, $X_n\subseteq X_{n+1}$ and $\bigcup_{n=1}^{\infty} X_n = X_{\infty}$) such that $\Sigma_n(X_n)=X_n$, we have $|X_n^t / \Sigma_n| \ge |X_{\infty}^t / \Sigma_{\infty}|$ for every large enough $n$.

However, not every exhaustion $\{X_n\}_{n=1}^{\infty}$ gives an asymptotically oligomorphic sequence: take for example $\Sigma_{\infty} \defeq S_{\infty}$ the group of finitely-supported permutations of $X_{\infty}\defeq\N$. 
Then $|X_{\infty}^t / \Sigma_{\infty}| < \infty$, but the exhaustions $\Sigma_n=S_n$ and $X_n = [2n]$ give $|X_n^t / \Sigma_n| \underset{n\to\infty}{\to} \infty$.
\end{remark}

\subsection{Word Measures on General Wreath Products}

The next result is a bound on the $w$-measure of the character $\textup{Ind}_{X_n}\phi$ for sequences of wreath products $G\wr_{X_n} \Sigma_n$, where now instead of $S_n\acts [n]$ we take general asymptotically oligomorphic group actions $\Sigma_n\acts X_n$.
The motivation for this result stems from the concept of stability: 

Recall from equation~\eqref{eq_main_res_of_ShoII} the informal discussion about stable characters, which lie at the heart of Conjecture~\ref{conj_HP22}.
We informally give some motivating facts about stable characters - see \cite{Sho23II} for the exact definitions and proofs.
For every sequence of groups $(\Sigma_n)_{n=1}^{\infty}$ and stable characters $\phi = (\phi_n)_{n=1}^{\infty}, \psi = (\psi_n)_{n=1}^{\infty}$ (where $\phi_n, \psi_n\in \textup{char}(\Sigma_n)$), the sequence of inner products stabilizes: $\inner*{\phi_n, \psi_n}$ is constant for large enough $n$.
Moreover, the trivial characters $(\textbf{1})_{n=1}^{\infty}$ form a stable (trivial) character.

Every sequence of actions $\Sigma_n\acts X_n$ yields a sequence of permutation characters $\psi_n$.
If this sequence is a stable character, then
$ \sup_n \left| X_n / \Sigma_n \right| = \sup_n \inner*{\psi_n, \textbf{1}} < \infty. $
This shows that stability is a stronger condition than being asymptotically oligomorphic.
Moreover, let $G$ be a compact group and let $\phi\in \hat{G}-\{\textbf{1}\}$. 
If $\Sigma_n\acts X_n$ yields a stable permutation character, then $\left( \textup{Ind}_{X_n} \phi \right)_{n=1}^{\infty}$ is a stable character of $\left(G\wr_{X_n} \Sigma_n\right)_{n=1}^{\infty}$ (see \cite{moeller2022extensions}). The word measure of such stable characters is controlled by the following theorem:

\begin{theorem}
\label{thm_general_wreath_character_bound}
Let $G$ be a compact group, let $\phi\in \hat{G}-\{\textbf{1}\}$, and let $\{\Sigma_n\acts X_n\}_{n=1}^{\infty}$ be an asymptotically oligomorphic sequence.
Then for every non-power word $w\in F_r$, the $w$-measure of $\textup{Ind}_n\phi\in \widehat{G\wr_{X_n} \Sigma_n}$ has polynomial decay rate in $\dim\left(\textup{Ind}_{X_n} \phi\right)$:
\[ \EX_w\left[\textup{Ind}_n\phi\right] = O\left(\dim(\textup{Ind}_n \phi)^{-1/2}\right). \]
\end{theorem}

\begin{conjecture}
\label{conj_reiter_for_word_measure_on_ao_grps}
    The $1/2$ above can be replaced by $1$ (which is tight).
\end{conjecture}

\noindent The proof combines the induction-convolution lemma with the results of \cite{reiter19}.
Recall the definition of \textbf{valid functions} from the proof of Lemma~\ref{lemma_ICL} (which was defined also in \cite[Definition 2.5: Associations]{reiter19}):

\begin{definition}
\label{def_valid_func}
 Let $\Sigma\acts X$ be a group action, $B\subseteq F_r$ a basis, $\Gamma = \Gamma_B(H)$ a core graph  of some f.g.\ subgroup $1\neq H\le F_r$ and let $\alpha\in \textup{Hom}(F_r, \Sigma)$.
 Then a function $f\colon V(\Gamma) \to X$ is called \textbf{valid} if for every $b$-labeled edge $\left(u_1\overset{b}{\to} u_2\right)\in E(\Gamma)$ we have $\alpha(b)(f(u_1)) = f(u_2)$. 
 
\noindent Informally, the following diagram commutes:
 \[\begin{tikzcd}
	{u_1} & {u_2} \\
	{f(u_1)} & {f(u_2)}
	\arrow["b", from=1-1, to=1-2]
	\arrow["f"', from=1-1, to=2-1]
	\arrow["f", from=1-2, to=2-2]
	\arrow["{\alpha(b)}", from=2-1, to=2-2]
\end{tikzcd}\]
\end{definition}

\noindent This definition is important due to the identity
\[\#\{x\in X \mid\,\, \forall h\in H: \alpha(h).x = x\} = \#\{\textup{valid functions }f\colon V(\Gamma) \to X\}. \]
\noindent This identity is easy to prove, and also follows from the proof of the induction-convolution lemma with $\phi=\textbf{1}$.

In \cite{reiter19}, Reiter bounded the probability that there exists a valid function on $\Gamma$, where $\Sigma\acts X$ is a finite, transitive group action and $\alpha\sim U(\textup{Hom}(F_r, \Sigma))$ is a uniformly random homomorphism. His proof uses the union bound (first moment method), thus giving a bound on the expected number of valid functions:
\begin{proposition}
(\cite[Corollary 5.6]{reiter19})
Let $\Sigma\acts X$ be a finite, transitive group action. 
Then for every $H\le F_r = \inner*{B}$ finitely generated and non-abelian with a core graph $\Gamma = \Gamma_B(H)$, 
and for every orbit $\pi \subseteq X^{V(\Gamma)}$ of the diagonal action $\Sigma\acts X^{V(\Gamma)}$, the expected number of valid functions in $\pi$ is at most $|X|^{-1/2}$.
\end{proposition}

Reiter also conjectured that the $\frac{1}{2}$ can be improved to $1$ (which is tight), and proved it for some special cases.

\begin{proof}
[Proof of Theorem~\ref{thm_general_wreath_character_bound}]
Let $\textup{orb}\defeq \sup_n \left| X_n^{|w|} / \Sigma_n \right|$ be the maximal number of orbits of functions $\Gamma_B(w)\to X_n$ (which is finite since $\Sigma_n\acts X_n$ is asymptotically oligomorphic).
For every $t\in \N$, the diagonal group action $\Sigma_n\acts X_n^t$ preserves the property of being injective:
the set of injective functions $[t]\hookrightarrow X_n$ is a union of $\Sigma_n$-orbits.
Denote the number of such orbits by $\textup{inj-orb}(t)$.
By the induction-convolution Lemma~\ref{lemma_ICL},
\begin{equation}
    \label{eq_quoting_ICL_yet_again}
    \EX_w\left[\textup{Ind}_n\phi\right] = \sum_{J\in \mathcal{Q}_B(w)} \EX_{w\to J}[\phi]\cdot L_{J\to F_r}^B(\Sigma_n\acts X_n).
\end{equation}
Recall the identity
\[L_{J\to F_r}^B(\Sigma_n\acts X_n) = \EX\left[\#\{\textup{valid injective functions }f\colon V(\Gamma_B(J)) \hookrightarrow X_n\}\right]. \]
Since $\EX_{w\to w}[\phi] = \inner*{\phi, 1} = 0$, and since $w$ is not a proper power, every contribution to equation \eqref{eq_quoting_ICL_yet_again} comes from a non-abelian (and of course finitely-generated) $J\in \mathcal{Q}_B(w)$, hence we may apply Reiter's theorem:
\[ L_{J\to F_r}^B(\Sigma_n\acts X_n)
\le \frac{\textup{inj-orb}\left(\left|V(\Gamma_B(J))\right|\right)}{\sqrt{|X_n|}}. \]
Moreover, every valid function $f\colon V\left(\Gamma_B(w)\right)\to X_n$ corresponds to a unique pair $\left(J, \Tilde{f}\right)$ where $J\in \mathcal{Q}_B(w)$ and $\Tilde{f}\colon V\left(\Gamma_B(J)\right)\hookrightarrow X_n$ is an injective valid function (this is how we built $\mathcal{Q}_B(w)$ in the first place), so
\[ \textup{orb} = \sum_{J\in \mathcal{Q}_B(w)} \textup{inj-orb}\left(\left|V\left(\Gamma_B(J)\right)\right|\right). \]
Using the naive bound $\left|\EX_{w\to J}[\phi]\right| \le \dim(\phi)$, we deduce 
\begin{equation*}
    \begin{split}
        \left|\EX_w\left[\textup{Ind}_n\phi\right] \right|
        &\le \dim(\phi)\cdot \sum_{J\in \mathcal{Q}_B(w), \, J\neq \inner{w}} L_{J\to F_r}^B(\Sigma_n\acts X_n)\\
        &\le \frac{\dim(\phi)}{\sqrt{|X_n|}}\cdot \sum_{J\in \mathcal{Q}_B(w), \, J\neq \inner{w}} \textup{inj-orb}\left(\left|V(\Gamma_B(J))\right|\right)\\
        &\le \frac{\dim(\phi)}{\sqrt{|X_n|}}\cdot \textup{orb}.
    \end{split}
\end{equation*}
Recalling that $\dim(\textup{Ind}_n \phi) = |X_n|\cdot \dim(\phi)$, we are done.
\end{proof}

\section{Free Product of Cyclic Groups}
\label{appendix_free_product_of_cyclic}

By following the proof of the induction-convolution lemma (Lemma~\ref{lemma_ICL}), one sees that the dependence on the distribution of $\sigma_1, \ldots, \sigma_r\in \Sigma_n $ is weak. 
For a probability measure $\mathcal{D}\colon \Sigma_n \to [0, 1]$, the product measure $ U\times\mathcal{D}\colon G^n\rtimes \Sigma_n \to [0, 1]$ (where $U$ is the Haar (uniform) measure on $G^n$) still satisfies an induction-convolution lemma: if $L_{\Gamma}^{\mathcal{D}, B}(n)$ is the expectation with respect to $((v_1, \sigma_1), \ldots, (v_r, \sigma_r))\sim (U\times\mathcal{D})^r$ of the number of valid injective functions $i\colon V(\Gamma)\hookrightarrow [n]$, then 
\[ \EX_{\left(v_j, \sigma_j\right)_{j=1}^r\sim \left(U\times \mathcal{D}\right)^r}\left[\textup{Ind}_n\phi\left(w\left(v_j, \sigma_j\right)_{j=1}^r\right)\right] = \sum_{\Gamma\in \mathcal{Q}_B(w)} \EX_{w\to H}[\phi]\cdot L_{\Gamma}^{B, \mathcal{D}}(n). \]
For example, $\mathcal{D}$ could be the uniform distribution on derangements in $S_n$.

This generalization is useful for the following problem: let $m\in \N$ and consider the free product $\mathcal{F}\defeq \inner{b_1, \ldots, b_r \divides b_1^m = \ldots = b_r^m = 1} = C_m * \ldots * C_m$.
Draw a uniformly random homomorphism $\alpha\sim U(\textup{Hom}(\mathcal{F}, S_n))$, and then for every $\gamma\in \mathcal{F}$ consider the "word" measure $\mu_{\gamma}\colon S_n\to [0, 1]$ obtained from $\mu_{\gamma}(\sigma) = \PR_{\alpha}(\alpha(\gamma)=\sigma)$. What is the expected number of fixed points with respect to this measure? according to \cite[Theorem 1.4]{PZimhony}, if $\gamma\in \mathcal{F}$ has infinite order,
\[ \lim_{n\to \infty} \EX_{\alpha}\left[\#\textup{fix}_n(\alpha(\gamma))\right] = \#\{\textup{subgroups }\mathcal{H}\le \mathcal{F} \textup{ of Euler characteristic }0\}. \]

In \cite[Conjecture 7.1]{PZimhony}, the authors defined for every $\gamma\in \mathcal{F}$
\[ \chi^{\max}(\gamma)\defeq \max\left\{ \chi(H) \mid \gamma\in H\le \mathcal{F}, \,\,\,\inner{\gamma} \textup{ is not a free factor isomorphic to }\Z \textup{ of }H \right\}, \]
with $\textup{Crit}(\gamma)$ the set of subgroups achieving the maximum, and conjectured that
\[\EX_{\alpha}\left[\#\textup{fix}_n(\alpha(\gamma))\right] = 1 + |\textup{Crit}(\gamma)|\cdot n^{\chi^{\max}(\gamma)} \left(1 + O\left(n^{-1/m}\right)\right). \] 

As in the free group scenario, the wreath product analogue turns out to be easier to solve:
\begin{definition}
We define
$\chi_{\mathcal{F}, \phi}(\gamma) \defeq \max\{\chi(\mathcal{H}): \gamma\in \mathcal{H}\le \mathcal{F}, \,\,\EX_{\gamma\to \mathcal{H}}[\phi]\neq 0 \}. $
Moreover, denote by $\textup{Crit}_{\mathcal{F}, \phi}$ the set of subgroups that achieve the maximum, and define 
\[ \mathscr{C}_{\mathcal{F}, \phi}(\gamma)\defeq \sum_{\mathcal{H}\in \textup{Crit}_{\mathcal{F}, \phi}} \EX_{\gamma\to \mathcal{H}}[\phi].\]    
\end{definition}

As in Proposition~\ref{prop_equiv_def_pi_phi}, for every $\gamma\in \mathcal{F}$ we have $\chi_{\mathcal{F}, \phi}(\gamma)\le \chi^{\max}(\gamma)$.

\begin{theorem}
Let $G$ be a finite group with $\textup{gcd}(|G|, m)=1$ (so that the map $x\mapsto x^m$ is invertible on $G$), and let $\phi\in \hat{G}-\{\textbf{1}\}$.
Then for every $\gamma\in \mathcal{F}$, the expectation of $\textup{Ind}_n\phi$ with respect to the measure of $\alpha(\gamma)$ where $\alpha\sim U(\textup{Hom}(\mathcal{F}, G\wr S_n))$ is 
\[ \EX_{\gamma}[\textup{Ind}_n\phi] = n^{\chi_{\mathcal{F}, \phi}(\gamma)}\cdot \left(\mathscr{C}_{\mathcal{F}, \phi}(\gamma) + O\left(n^{-1/m}\right)\right). \]
\end{theorem}

This is the analogue of \cite[Conjecture 7.1]{PZimhony} stated for wreath products.

\begin{proof}
Let $\mathcal{D} = \mu_{b_1}$ be the distribution on $S_n$ that gives a random solution of the equation $X^m = 1$. 
Let $\gamma\in \mathcal{F}$. 
Choose a representative $w\in F_r$ which is sent to $\gamma$ under the quotient map $F_r\to F_r/\inner{\inner{b_1^m, \ldots, b_r^m}} = \mathcal{F}$ such that all of the exponents of letters in $w$ are between $0, \ldots, m-1$.
Now the $\gamma$-measure on $G\wr S_n$ is the pushforward of the measure $(U\times \mathcal{D})^r$ on $(G\wr S_n)^r$ with respect to the word map $w\colon (G\wr S_n)^r\to G\wr S_n$, and we can apply the induction-convolution lemma: 
\[ \EX_{\gamma}[\textup{Ind}_n\phi] = \EX_{\left(v_j, \sigma_j\right)_{j=1}^r\sim \left(U\times \mathcal{D}\right)^r}\left[\textup{Ind}_n\phi\left(w\left(v_j, \sigma_j\right)_{j=1}^r\right)\right] = \sum_{\Gamma\in \mathcal{Q}_B(w)} \EX_{w\to H}[\phi]\cdot L_{\Gamma}^{B, \mathcal{D}}(n). \]
By \cite[Theorem 2.6]{PZimhony}, for every $\Gamma\in \mathcal{Q}_B(w)$, either $ L_{\Gamma}^{B, \mathcal{D}}(n) = 0$ or
\[ L_{\Gamma}^{B, \mathcal{D}}(n) = n^{\chi(\Gamma)}\cdot \left(1 + O\left(n^{-1/m}\right)\right), \]
and the result follows as in the proof of Theorem~\ref{thm_main_result}.
\end{proof}

Since the formulation and notation in \cite{PZimhony} are very different, we explain now how to derive this conclusion from \cite[Theorem 2.6]{PZimhony}.
The formulation in \cite{PZimhony} is as follows:
fix a presentation complex $X_{\mathcal{F}}$ of $\mathcal{F}$. 
For every compact sub-cover $p\colon Y\to X_{\mathcal{F}}$ (that is, a compact CW sub-complex of a finite cover of $X_{\mathcal{F}}$), denote by $\EX^{\textup{emb}}_Y(n)$ the expected number of injective lifts of $Y$ to a random covering of $X_{\mathcal{F}}$ of degree $n$. Then
\[ \EX^{\textup{emb}}_Y(n) = n^{\chi(p_*(\pi_1(Y)))}\cdot \left(1 + O\left(n^{-1/m}\right)\right). \]

Similarly to the poset $\mathcal{Q}_B(w)$ of quotients of $\Gamma_B(w)$, there is a poset $\mathcal{R}_{\gamma}$ of all surjective lifts of $\gamma\colon \Ss^1\twoheadrightarrow Y$ to sub-covers of $ X_{\mathcal{F}}$.
As in the free group setting, where $\EX_w[\#\textup{fix}_n] = \sum_{H\in \mathcal{Q}_B(w)} L_H^B(n)$, we have
\[ \EX_{\alpha}[\#\textup{fix}(\alpha(\gamma))] = \sum_{Y\in \mathcal{R}_{\gamma}} \EX^{\textup{emb}}_Y(n). \]

The case of $C_m * \ldots * C_m$ is simpler than the more general free product handled by \cite{PZimhony}.
Fix $\gamma, w$ as in the proof above.
Then there is a bijection between $\mathcal{R}_{\gamma}$ and the set of graphs $\Gamma\in \mathcal{Q}_B(w)$ such that $L_{\Gamma}^{B, \mathcal{D}}(n)\neq 0$, and for every such graph $\Gamma$ we have $\EX^{\textup{emb}}_Y(n) = L_{\Gamma}^{B, \mathcal{D}}(n)$ where $Y\in \mathcal{R}_{\gamma}$ is the corresponding sub-cover.

\section{Open Questions}
\label{section_open_questions}

We conclude with some open questions naturally arising from the results in this paper.

\begin{enumerate}
\item Let $G$ be a compact group and $\phi\in \hat{G}-\{\textbf{1}\}$. Let $w\in H\le J\le F_r$ be a chain of $B$-surjective extensions (for some basis $B\subseteq F_r$).
By Remark~\ref{remark_composition_with_balanced_morphisms}, if $\phi$ is linear then 
\begin{itemize}
    \item If $\EX_{\eta_{w\to H}}[\phi]\neq 0$, then also $\EX_{\eta_{w\to J}}[\phi]\neq 0$. 
    \item For every $w\in F_r$, $\pi_{\phi}(w)\in \{0, 1, \ldots, r\} \sqcup \{\infty\}$. (This follows from the previous item).
\end{itemize}
We have also seen in Remark~\ref{remark_composition_with_balanced_morphisms} that for every $n\ge 3$, the first property does not hold for the standard character of $S_n$.
It is natural to ask: what are the possible values of $\pi_{\phi}(w)$ when $\phi$ is not linear? 
Which $\phi$ satisfies $\pi_{\phi}(w)\in \{0, 1, \ldots, r\} \sqcup \{\infty\}$ for every $w$?
\item 
    (Reiter's Conjecture) Let $H\le F_r$ be a finitely-generated non-abelian subgroup, and let $(\Sigma_n\acts X_n)_{n=1}^{\infty}$ be an asymptotically oligomorphic sequence of group actions. 
    Then
     $ \EX_{\eta_{H\to F_r}}[\Sigma_n\acts X_n] = O\left(|X_n|^{-1}\right). $
\end{enumerate}

\section{Glossary}
\begin{table}[ht!]
\centering
\begin{tabular}{||c | c | c ||} 
 \hline
 Notation 
 & Description 
 & Comments \\[0.7ex] 
 \hline\hline
 $[n]$ & $\{1, \ldots, n\}$ & \\[0.7ex] 
 $(n)_t$ & $ n\cdot (n-1)\cdots (n-t+1)$ & Falling factorial \\[0.7ex]
 $F_r$ & $\textup{Free}(\{b_1, \ldots, b_r\})$ & The ambient free group \\[0.7ex]
 $\Omega_r$ & $\Gamma_{\{b_1, \ldots, b_r\}}(F_r) $ & The bouquet with $r$ petals \\[0.7ex]
 $\eta$ & Morphism of core graphs & Also inclusion of free subgroups\\[0.7ex]
 $\hat{G}$ & Irreducible characters of $G$ &  \\[0.7ex] 
 $\textup{char}(G)$ & All characters of $G$ &  \\[0.7ex]
 $\textup{conj}(G)$ & Conjugacy classes of $G$ &\\[0.7ex]
 $C_m$ & $\Z/m\Z$ & $C_{\infty} = \Ss^1$\\[0.7ex]
 $\textbf{1}$ & $\textbf{1}\colon G\to \{1\}$ & The trivial character\\[0.7ex]
 $w$ & Element of $F_r$&\\[0.7ex]
 $|w|_B$ & The length of $w$ in basis $B$ & \\[0.7ex]
 $|w|$ & $\min_B |w|_B$ & \\[0.7ex]
 $\EX_w[f]$ & $\EX_{\alpha\sim U(\textup{Hom}(F_r, G))} [f(\alpha(w))]$ & The $w$-measure of $f$ \\ [1ex] 
 \hline
\end{tabular}
\caption{Glossary}
\label{table:glossary}
\end{table}
\FloatBarrier

\printbibliography 
\noindent
Yotam Shomroni,\\
School of Mathematical Sciences,\\
Tel Aviv University,\\
Tel Aviv, 6997801, Israel\\
yotam.shomroni@gmail.com

\end{document}